\documentclass{amsart}

\usepackage[a4paper,hmargin=2.8cm,vmargin=4cm]{geometry}
\usepackage{amsfonts,amssymb,amscd,amstext}
\usepackage{mathtools}
\mathtoolsset{showonlyrefs=true} 
\usepackage[utf8]{inputenc}
\usepackage[colorlinks=true,linkcolor=blue,citecolor=red]{hyperref}
\usepackage{braket,relsize,nccmath}

\makeatletter
\@namedef{subjclassname@2020}{%
  \textup{2020} Mathematics Subject Classification}
\makeatother

\renewcommand{\leq}{\leqslant}
\renewcommand{\geq}{\geqslant}
\newcommand{\ptl}{\partial}
\newcommand{\rr}{{\mathbb{R}}}
\newcommand{\rrn}{\mathbb{R}^{n+1}}
\newcommand{\la}{\lambda}
\newcommand{\sph}{\mathbb{S}}
\newcommand{\nn}{\mathbb{N}}
\newcommand{\sub}{\subset}
\newcommand{\escpr}[1]{\big<#1\big>}
\newcommand{\Sg}{\Sigma} 
\newcommand{\Om}{\Omega}
\newcommand{\eps}{\varepsilon}
\newcommand{\var}{\varphi}
\newcommand{\ind}{\mathcal{I}}
\newcommand{\indo}{\mathcal{I}}
\newcommand{\cc}{\mathcal{C}}
\newcommand{\hh}{\mathcal{H}}
\newcommand{\ovn}{\overline{N}}
\newcommand{\ar}{\mathcal{A}}
\newcommand{\bb}{\mathcal{B}}
\newcommand{\pp}{\mathcal{P}}

\DeclareMathOperator{\divv}{div}

\newtheorem{theorem}{Theorem}[section]
\newtheorem{proposition}[theorem]{Proposition}
\newtheorem{lemma}[theorem]{Lemma}
\newtheorem{corollary}[theorem]{Corollary}

\theoremstyle{definition}

\newtheorem{remark}[theorem]{Remark}
\newtheorem{remarks}[theorem]{Remarks}
\newtheorem{example}[theorem]{Example}

\theoremstyle{remark}

\numberwithin{equation}{section}

\AtBeginDocument{%
   \def\MR#1{}
}

\begin{document}

\title[Anisotropic partitioning problem in convex domains]{On the anisotropic partitioning problem \\ in Euclidean convex domains}

\author[C\'esar Rosales]{C\'esar Rosales}
\address{Departamento de Geometr\'{\i}a y Topolog\'{\i}a and Excellence Research Unit ``Modeling Nature'' (MNat) Universidad de Granada, E-18071,
Spain.} 
\email{crosales@ugr.es}

\date{February 25, 2025}

\thanks{This research was supported by the grant PID2020-118180GB-I00 funded by MCIN/AEI/10.13039/501100011033 and by the grant PID2023-151060NB-I00 funded by MCIN/AEI/10.13039/501100011033 and by ERDF/EU} 



\begin{abstract}
We consider the variational problem of minimizing an anisotropic perimeter functional under a volume constraint in a Euclidean convex domain. We extend to this setting analytical properties of the isoperimetric profile, topological features about the minimizers and sharp isoperimetric inequalities with respect to convex cones. Besides some geometric measure theory results about the existence and regularity of minimizers, the proofs rely on a second variation formula for the anisotropic area of hypersurfaces with non-empty boundary.
\end{abstract}

\maketitle

\thispagestyle{empty}

\section{Introduction}
\label{sec:intro}
\setcounter{equation}{0}

The \emph{anisotropic isoperimetric problem} (or \emph{Wulff problem}) in a Euclidean domain seeks those sets minimizing a perimeter functional associated to a surface tension while enclosing a fixed volume. The solutions are called \emph{anisotropic isoperimetric regions} or \emph{minimizers}. Some related questions are the existence of minimizers, the regularity of their boundaries and the analytical properties of the corresponding isoperimetric profile. The obtention of geometric, topological and variational information is also interesting because, in general, the complete description of the solutions is a challenging task. 

In this work, for a convex body $K\subset\rrn$ with $0\in K^\circ$, and a domain $\Om\subset\rrn$, we consider the \emph{anisotropic perimeter} $\pp_K(\cdot\,,\Om)$ obtained when the support function $h_K$ is used as an elliptic integrand. More precisely, for any set $E$ of locally finite perimeter in $\Om$, we have
\[
\pp_K(E,\Om)=\int_{\ptl^*E}h_K(N_E)\,dA,
\]
where $\ptl^*E$ is the reduced boundary in $\Om$, $N_E$ is the measure-theoretic outer unit normal, and $dA$ is the Euclidean area element. We note that $\ptl E\cap\ptl\Om$ has no contribution to $\pp_K(E,\Om)$. When $E$ is open and $\ptl E\cap\Om$ is a smooth hypersurface $\Sg$, then $\pp_K(E,\Om)$ coincides with the \emph{anisotropic area}
\[
\ar_K(\Sg):=\int_{\Sg}h_K(N)\,dA,
\]
where $N$ is the outer normal on $\Sg$. In the isotropic case, where $K$ is the round unit ball about $0$, we get that $h_K$ is the Euclidean norm, and we recover the notion of relative perimeter. For an arbitrary $K$ as above, the function $h_K$ provides an \emph{asymmetric norm} which may fail to be a norm. Thus, the anisotropic area of a two-sided hypersurface $\Sg$ may depend on the chosen unit normal $N$. 

Since the anisotropic isoperimetric problem is a non-trivial generalization of the classical one, it is natural to extend important results in the isotropic setting by assuming different geometric and analytical conditions on $K$ and $\Om$. Though we will focus on proper domains, it is worth mentioning that the question was completely solved in $\rrn$ by Taylor~\cite{taylor-unique,taylor}, see also Fonseca and Müller~\cite{fonseca-muller}, and Brothers and Morgan~\cite{brothers-morgan}. They showed that, up to translations, dilations about $0$, and sets of volume zero, the convex body $K$ uniquely minimizes the perimeter $\pp_K(\cdot\,,\rrn)$ among sets of the same volume. The boundary $\ptl K$ is usually referred to as the \emph{Wulff shape} in honor of the crystallographer Georg Wulff, who first constructed the optimal crystal for a specified integrand~\cite{wulff}. After previous work of Morgan~\cite{morgan-alexandrov} for arbitrary planar norms, the characterization of the critical sets in $\rrn$ via an Alexandrov type theorem was established by He, Li, Ma and Ge~\cite{alexandrov-anisotropic} for smooth sets, and by De Rosa, Kolasi\'nksi and Santilli~\cite{rks} for finite perimeter sets. On the other hand, the uniqueness of the Wulff shape as a compact, immersed, two sided, anisotropic stable hypersurface (second order minimum of $\ar_K$ for fixed volume) was given by Palmer~\cite{palmer} and Winklmann~\cite{winklmann}.

In a half-space $\hh$, Winterbottom~\cite{winterbottom} proved that a suitable truncation of $K$ minimizes a capillarity functional involving not only the anisotropic perimeter but also the wetting area in $\ptl\hh$. The corresponding classifications of the smooth critical sets and of the compact stable hypersurfaces immersed in $\hh$ were treated by Jia, Wang, Xia and Zhang~\cite{jwxz} and by Guo and Xia~\cite{guo-xia}, respectively. We must mention a paper of Baer~\cite{baer}, who employed anisotropic symmetrization to deduce existence, convexity and symmetry of minimizers in $\hh$ when $K$ is a smooth, axially symmetric and strictly convex body. Under similar conditions, Koiso and Pamer~\cite{koiso-palmer2,koiso-palmer4,koiso-palmer3,kp-circular} analyzed the compact stable hypersurfaces for the anisotropic capillarity problem in a slab $\mathcal{S}\subset\rr^3$. In a convex cone $\cc\subset\rrn$ with vertex at the origin, Milman and Rotem~\cite{milman-rotem}, and Cabr\'e, Ros-Oton and Serra~\cite{cabre-rosoton-serra}, proved that $\lambda K\cap\cc$ is an anisotropic isoperimetric region for any $\lambda>0$. The characterization of the minimizers in $\cc$ was later provided by Dipierro, Poggesi and Valdinoci \cite{dipierro-poggesi-valdinoci}. The smooth critical sets were discussed by Weng~\cite{weng}, and by Jia, Wang, Xia and Zhang~\cite{jwxz2}. The paper \cite{jwxz2} contains a Heintze-Karcher inequality in arbitrary convex domains which, in the case of a convex cone $\cc$, leads to an Alexandrov theorem for smooth critical sets avoiding the singular part $(\ptl\cc)_S$ of $\ptl\cc$. As to the classification of compact stable hypersurfaces avoiding $(\ptl\cc)_S$, we refer the reader to~\cite{rosales-cones} after previous work of Koiso~\cite{koiso-equilibrium,koiso-2023} for wedge-shaped domains. We also refer to Barbosa and Carvalho Silva~\cite{barbosa-silva}, who investigated critical and stable anisotropic hypersurfaces with respect to some axially symmetric surface tensions in smooth domains of revolution in $\rr^3$.

In this article we consider arbitrary convex domains and derive anisotropic counterparts to relevant results concerning analytical properties of the isoperimetric profile, topological features about the minimizers and sharp isoperimetric comparisons with respect to convex cones. Let us describe them in more detail.

For a bounded domain $\Om\subset\rrn$, the \emph{anisotropic isoperimetric profile} $I_{\Om,K}$ is the function which provides the least possible perimeter $\pp_K(E,\Om)$ among sets $E\subset\Om$ of fixed volume, see equation~\eqref{eq:profile}. Basic statements including the positivity and the continuity of $I_{\Om,K}$ are found in Section~\ref{sec:partprob}. In Theorem~\ref{th:concavity} we show the concavity of the function $\psi:=(I_{\Om,K})^{(n+1)/n}$ for any bounded convex domain $\Om\subset\rrn$. When $K$ is the round unit ball and $\Om$ is smooth, this was established by Kuwert~\cite{kuwert}, who improved the concavity of $I_{\Om,K}$ discovered by Sternberg and Zumbrun~\cite{sz}. The extension to non-smooth domains can be achieved by approximation arguments, see Milman~\cite[Sect.~6]{milman} and Leonardi, Ritor\'e and Vernadakis~\cite[Sect.~5]{unbounded}. Indeed, the concavity of $\psi$ still holds in the isotropic setting for unbounded convex domains of uniform geometry~\cite[Sect.~5]{unbounded}. The aforementioned approximation argument together with the existence of anisotropic minimizers in bounded domains allows to deduce Theorem~\ref{th:concavity} from Theorem~\ref{th:fundamental}, where we see that the upper second derivative of $\psi$ satisfies the inequality $\overline{D^2}\psi(v_0)\leq 0$ in any smooth convex domain $\Om$ having minimizers of volume $v_0$. Our tools permit to discuss the equality $\overline{D^2}\psi(v_0)=0$, which implies strong rigidity conditions on the minimizers of volume $v_0$. In Corollary~\ref{cor:unique1}, see also Remark~\ref{re:unique2}, we employ such conditions to recover uniqueness statements in $\rrn$ and smooth convex cones. 

From the concavity of $\psi$ it is straightforward to get some analytical consequences (Corollary~\ref{cor:regprop}), as the strict subbaditivity of $I_{\Om,K}$. This entails the connectivity of the anisotropic minimizers (Theorem~\ref{th:connectset}) by using the same reasoning as in the isotropic setting, see Bayle and the author~\cite{bayle-rosales}. After this, we analyze a stronger topological property: the connectivity of the interior boundary of any minimizer. In Theorem~\ref{th:connectbound} we prove that, in a smooth convex domain $\Om$, and with respect to a smooth strictly convex body $K$, the interior boundary of any anisotropic isoperimetric region is connected unless it satisfies some rigidity conditions. A similar statement in the isotropic case was given by Sternberg and Zumbrun~\cite{sz} when $\Om$ is bounded. We note that Theorem~\ref{th:connectbound} does not prevent connected minimizers with disconnected interior boundary. This motivates us to find in Corollary~\ref{cor:connected24} additional hypotheses ensuring the connectivity of the isoperimetric interior boundaries.

Finally, we obtain an anisotropic isoperimetric inequality in a bounded convex domain $\Om$. More precisely, we establish in Theorem~\ref{th:maincomp} the profile comparison $I_{\Om,K}\leq I_{\cc,K}$, where $\cc$ is any open convex cone containing $\Om$ and with vertex at $\ptl\Om$. This estimate is suggested by the concavity of the function $\psi=(I_{\Om,K})^{(n+1)/n}$ and the linearity of $(I_{\cc,K})^{(n+1)/n}$, see equation~\eqref{eq:isoprofcones}. According to this fact, we only have to check the desired comparison for small volumes. And this is done in Proposition~\ref{prop:1stcomp} for arbitrary domains under mild convexity assumptions. We observe that the inequality $I_{\Om,K}\leq I_{\cc,K}$ is sharp, in the sense that, if equality holds for some positive volume, then $\ptl\Om$ and $\ptl\cc$ locally coincide around the vertex of $\cc$. In the isotropic case the theorem was proved by Kuwert~\cite{kuwert} for smooth domains. The extension to non-smooth ones is due to Ritor\'e and Vernadakis~\cite[Sect.~6]{bounded} when $\Om$ is bounded, see Leonardi, Ritor\'e and Vernadakis~\cite[Sect.~6.3]{unbounded} for the unbounded case.

The proofs of our results are inspired in the isotropic setting. The main ingredients are anisotropic versions of standard tools in geometric measure theory and a second variation formula for the anisotropic area $\ar_K$. On the one hand, the existence and regularity of solutions to the anisotropic isoperimetric problem are discussed in Section~\ref{sec:partprob} by giving precise references. In a bounded Lipschitz domain $\Om\subset\rrn$ there are anisotropic minimizers of any volume. Moreover, if $\Om$ is smooth, then the interior boundary of any minimizer enjoys the same analytical properties as in the isotropic case, with the difference that the Hausdorff dimension of the singular set is at most $n-2$. Though this estimate is worse than the classical one, it is enough to prove Proposition~\ref{prop:approx}; this provides a useful approximation argument for describing the anisotropic minimizers by focusing on the regular part of the interior boundary. On the other hand, in Section~\ref{sec:2ndvar} we compute $\ar_K''(0)$ for a two sided smooth hypersurface $\Sg$ (with non-empty boundary) when this is deformed under a flow having velocity vector field $X$ proportional to the \emph{anisotropic normal} $N_K$. When $\Sg$ is the regular part of the interior boundary of an anisotropic minimizer, the equilibrium conditions described in Proposition~\ref{prop:varprop} lead to simplify $\ar_K''(0)$ for volume preserving variations. The resulting expression in Proposition~\ref{prop:2ndareavol} generalizes the classical one by means of a quadratic form $\indo_K$ involving some of the anisotropic analytical and geometric notions introduced in Section~\ref{subsec:extrinsic}. The study of $\indo_K$ should be helpful in the analysis of the associated stability problems. Previous second variational formulas in anisotropic capillarity problems were found by Koiso and Palmer~\cite{koiso-palmer2,koiso-palmer4}, De Philippis and Maggi~\cite{philippis-maggi-2}, Barbosa and Carvalho Silva~\cite{barbosa-silva}, and Guo and Xia~\cite{guo-xia}. Though our formula is only valid for compactly supported deformations of the regular set $\Sg$, we will be able to infer geometric information on the minimizers via the approximation argument in Proposition~\ref{prop:approx}.

Our results apply in the particular situation where $K$ is a centrally symmetric convex body. In this case we gain information from the fact that the support function $h_K$ defines a norm. For instance, the anisotropic area $\ar_K$ does not depend on the chosen normal, the isoperimetric profile $I_{\Om,K}$ is symmetric with respect to $V(\Om)/2$, and the complementary set $\Om\setminus E$ of a minimizer is a new minimizer. More consequences are shown in Corollaries~\ref{cor:regprop} and~\ref{cor:connected24}. It is also remarkable that the anisotropic profile $I_{\hh,K}$ is the same for any open half-space $\hh\subset\rrn$, which leads to improve the comparison in Theorem~\ref{th:maincomp}, see Example~\ref{ex:halfspace}. Of course, none of these properties necessarily hold within an arbitrary convex body $K$. We must observe that the centrally symmetric case is also covered in recent papers of Antonelli, Pasqualetto, Pozzetta and Semola~\cite{antonelli,antonelli2}. These authors considered the more general setting of $\text{RCD}(q,N)$ metric measure spaces $(X,d,\mathcal{H}^N)$ with a uniform lower bound on the volume of unit balls. They established the inequality $\psi''\leq 0$ in viscosity sense and deduced several consequences, including fine properties of the isoperimetric profile and connectedness of minimizers. Their approach rely on generalized existence of isoperimetric regions and a Laplacian comparison for the distance function from the boundary, which plays the same role as the first and second variation formulas for the perimeter of equidistant sets.

The techniques employed in the paper allow also to study \emph{anisotropic minimal hypersurfaces with free boundary} in a Euclidean domain $\Om$. These are the critical points of $\ar_K$ under compactly supported deformations leaving invariant $\ptl\Om$. When the hypersurface $\Sg$ is also \emph{stable}, i.e., it is a second order minimum of $\ar_K$ for such variations, then it satisfies the functional inequality $\ind_K(\omega)\geq 0$ for any $\omega\in C^\infty_0(\Sg)$. From here, and assuming that $\Om$ is smooth and convex, we can infer in Theorem~\ref{th:minimal} the same rigidity conditions as in Theorem~\ref{th:connectbound}. 

Besides this introduction,  we have organized the paper into four sections and an appendix. In Section~\ref{sec:prelimi} we gather some notions and facts concerning the relative geometry of hypersurfaces and the anisotropic perimeter of sets. In Section~\ref{sec:partprob} we introduce the anisotropic partitioning problem in Euclidean domains and review existence, regularity and approximation properties for the solutions.  Section~\ref{sec:2ndvar} is devoted to compute the second variation of the anisotropic area for certain deformations of a hypersurface with non-empty boundary. The main results of the paper are contained in Section~\ref{sec:main}. Finally, in Appendix~\ref{app:minimal} we discuss anisotropic minimal hypersurfaces with free boundary.

\section{Preliminaries}
\label{sec:prelimi}
\setcounter{equation}{0}

In this section we collect definitions and results that we need throughout this work. We will follow the point of view and notation in \cite[Sect.~2]{rosales-cones}. A more general approach in the setting of relative differential geometry is found in Reilly~\cite[Sect.~1]{reilly-anisotropic}.

\subsection{Some facts about convexity}
\label{subsec:convexity}
\noindent

By a \emph{convex body} (about the origin) we mean a compact convex set $K\subset\rrn$ such that $0\in K^\circ$. The associated \emph{support function} $h_K:\rrn\to\rr$ is given by
\begin{equation}
\label{eq:hK}
h_K(w):=\max\{\escpr{u,w}\,;\,u\in K\}, 
\end{equation}
where $\escpr{\cdot\,,\cdot}$ is the standard scalar product in $\rrn$. This provides an \emph{asymmetric norm}~\cite[Sect.~1.7.1]{schneider} which may fail to be a norm. In particular, the map $d_K(p,q):=h_K(p-q)$ need not be a distance in $\rrn$ because the symmetry property need not hold. In the special case where $K$ is centrally symmetric about $0$ we have that $h_K$ is a norm and $d_K$ is a distance. Anyway, the function $h_K$ is always convex and so, it is continuous on $\rrn$ and Lipschitz on compact sets~ \cite[Thm.~1.5.3]{schneider}. 

For any $w\neq 0$, the \emph{support set} of $K$ at $w$ is $\Pi_w\cap K$, where $\Pi_w:=\{p\in\rrn\,;\,\escpr{p,w}=h_K(w)\}$ is the \emph{supporting hyperplane} of $K$ with exterior normal $w$. If $\Pi_w\cap K$ is a single point $p$, then $h_K$ is differentiable at $w$, and its gradient satisfies $(\nabla h_K)(w)=p$, see \cite[Cor.~1.7.3]{schneider}. 

Suppose now that $K$ is a \emph{strictly convex body}, which means that the topological boundary $\ptl K$ does not contain segments. For any $w\neq 0$ there is a point $\pi_K(w)\in\ptl K$ for which $\Pi_w\cap K=\{\pi_K(w)\}$. This defines a map $\pi_K:\rrn_*\to\ptl K$, called the \emph{$K$-projection}, such that $\escpr{\pi_K(w),w}=h_K(w)$ for any $w\neq 0$. It follows that $h_K$ is differentiable in $\rrn_*$ and
\begin{equation}
\label{eq:gradhK}
(\nabla h_K)(w)=\pi_K(w), \quad\text{for any }  w\neq 0.
\end{equation}
If we further assume that $K$ is smooth, i.e., $\ptl K$ is a $C^\infty$ hypersurface, then the associated outer unit normal map $\eta_K:\ptl K\to\sph^n$ is a diffeomorphism and $\pi_K(w)=\eta_K^{-1}(w/|w|)$ for any $w\neq 0$ (here $|\cdot|$ stands for the Euclidean norm and $\sph^n$ for the set of vectors $w$ with $|w|=1$). Thus $h_K$ is smooth ($C^\infty$) on $\rrn_*$. On the other hand, for any $w\in\sph^n$ the equality $(\eta_K\circ\pi_K)(w)=w$ entails that the differential at $w$ of the map $\pi_{K}:\sph^n\to\ptl K$ is the endomorphism $(d\pi_K)_w:w^\bot\to w^\bot$ given by $(d\pi_K)_w(v)=(d\eta_K)^{-1}_{\pi_K(w)}(v)$. Since $\ptl K$ is an ovaloid, this endomorphism is positive definite and symmetric with respect to the scalar product restricted to $w^\bot$. By continuity, we can find constants $a,b>0$ depending only on $K$ such that
\begin{equation}
\label{eq:bound}
a\,|v|^2\leq\escpr{(d\pi_K)_w(v),v}\leq b\,|v|^2,\quad\text{for any }w\in\sph^n \text{ and } v\in w^\bot.
\end{equation}
We observe that, in the isotropic case (where $K$ is the round unit ball about $0$), we have $h_K(w)=|w|$ and $\pi_K(w)=w/|w|$ for any $w\neq 0$, whereas $\eta_K(p)=p$ for any $p\in\ptl K=\sph^n$.

Next, we recall some facts about the Hausdorff distance in $\rrn$ that will be useful for the approximation argument in the proof of Theorem~\ref{th:concavity}. For two compact non-empty sets $C,D\subset\rrn$ their \emph{Hausdorff distance} is the number
\[
\delta(C,D):=\min\{r\geq 0\,;\,C\subseteq D+r\bb, D\subseteq C+r\bb\},
\]
where $\bb$ is the closed unit ball, $+$ is the Minkowski sum of sets and $r\bb:=\{r\,p\,;\,p\in\bb\}$. A sequence of compact non-empty sets $\{C_i\}_{i\in\nn}$ in $\rrn$ converges in Hausdorff distance to a compact non-empty set $C\subset\rrn$ if $\{\delta(C_i,C)\}\to 0$ when $i\to\infty$. In the next lemma we gather some properties of this convergence when the involved sets are convex, see \cite[Sect.~1.8]{schneider}.

\begin{lemma}
\label{lem:hausdorff}
Let $\{C_i\}_{i\in\nn}$ be a sequence of compact convex non-empty sets in $\rrn$ converging in Hausdorff distance to a compact convex non-empty set $C$. Then, we have
\begin{itemize}
\item[(i)] $\{V(C_i)\}\to V(C)$ when $i\to\infty$, where $V$ denotes the Euclidean volume.
\item[(ii)] For any $p\in C$ there is a sequence $\{p_i\}_{i\in\nn}$ such that $p_i\in C_i$ for any $i\in\nn$ and $\{p_i\}\to p$. Moreover, if $\{p_i\}_{i\in\nn}$ is a sequence with $p_i\in C_i$ for any $i\in\nn$ and $\{p_i\}\to p$, then $p\in C$. Indeed, these two properties characterize the Hausdorff convergence of $\{C_i\}_{i\in\nn}$ to $C$.  
\item[(iii)] If $A$ is a compact set with $A\subset C^\circ$, then there is $i_0\in\nn$ such that $A\subset C_i^\circ$ for any $i\geq i_0$.
\item[(iv)] If $0\in C^\circ$, then there is $i_0\in\nn$ and sequences $\{\theta_i\}_{i\geq i_0}\subset(0,1]$, $\{\mu_i\}_{i\geq i_0}\subset[1,\infty)$ with $\{\theta_i\}\to 1$, $\{\mu_i\}\to 1$ and $\theta_i\,C\subseteq C_i\subseteq\mu_i\,C$ for any $i\geq i_0$.
\end{itemize}
\end{lemma}

\subsection{Anisotropic extrinsic geometry}
\label{subsec:extrinsic}
\noindent

In this subsection we denote by $K\subset\rrn$ a smooth strictly convex body and by $\Sg\subset\rrn$ a two-sided smooth hypersurface, possibly with non-empty boundary $\ptl\Sg$. 

For a fixed smooth unit normal vector field $N$ on $\Sg$, the associated \emph{shape operator} at a point $p\in\Sg$ is the symmetric endomorphism $B_p:T_p\Sg\to T_p\Sg$ defined by $B_p(w):=-D_wN$. Here $T_p\Sg$ is the tangent hyperplane of $\Sg$ at $p$, while $D$ stands for the flat Levi-Civita connection in $\rrn$. 

The \emph{anisotropic Gauss map} or \emph{anisotropic normal} on $\Sg$ is the map $N_K:\Sg\to\ptl K$ given by
\begin{equation}
\label{eq:normal} 
N_K:=\pi_K\circ N.
\end{equation} 
By setting
\[
\var_K:=h_K\circ N,
\] 
we infer from the definition of $\pi_K$ that
\begin{equation}
\label{eq:nkn}
\escpr{N_K,N}=\var_K\quad\text{on }  \Sg. 
\end{equation}
When $\ptl\Sg\neq\emptyset$ we introduce the \emph{anisotropic conormal} along $\ptl \Sg$ by equality
\begin{equation}
\label{eq:conormal}
\nu_K:=\var_K\,\nu-\escpr{N_K,\nu}\,N,
\end{equation}
where $\nu$ is the \emph{inner conormal} along $\ptl\Sg$. Note that $\nu_K$ is a normal vector to $\ptl\Sg$ with $\escpr{N_K,\nu_K}=0$.

Fix a point $p\in\Sg$. As indicated in Section~\ref{subsec:convexity}, the differential at $N(p)$ of $\pi_K:\sph^n\to\ptl K$ is a positive definite and symmetric endomorphism of $N(p)^\bot=T_p\Sg$. We will simply write
\begin{equation}
\label{eq:endoq}
Q_p:=(d\pi_K)_{N(p)}.
\end{equation}
The \emph{anisotropic gradient} of a function $f\in C^\infty(\Sg)$ is the tangent vector field
\begin{equation}
\label{eq:angrad}
(\nabla_\Sg^K f)(p):=Q_p((\nabla_\Sg f)(p)), \quad p\in\Sg,
\end{equation}
where $\nabla_\Sg$ is the gradient in $\Sg$ for the induced metric. The \emph{anisotropic shape operator} at $p$ is the endomorphism
\begin{equation}
\label{eq:shape}
(B_K)_p:=-(dN_K)_p=Q_p\circ B_p.
\end{equation}
The \emph{anisotropic mean curvature} at $p$ is the number
\begin{equation}
\label{eq:mc}
H_K(p):=\frac{\text{tr}((B_K)_p)}{n}=-\frac{(\divv_\Sg N_K)(p)}{n},
\end{equation}
where $\text{tr}(\cdot)$ is the trace of endomorphisms in $T_p\Sg$ and $\divv_\Sg$ is the divergence relative to $\Sg$. When $K$ is the round unit ball about $0$ then $Q_p$ is the identity map in $T_p\Sg$, so that $(B_K)_p=B_p$ and $H_K(p)$ equals the Euclidean mean curvature of $\Sg$ at $p$.

The following statement is a useful characterization result, which generalizes a well known fact in the isotropic case, see Palmer~\cite[p.~3666]{palmer}, Clarenz~\cite[Proof of Thm.~4.1]{clarenz-2004} or \cite[Sect.~2]{rosales-cones}.

\begin{proposition}
\label{prop:umbilical}
Let $\Sg\subset\rrn$ be a two-sided smooth hypersurface. Then, we have
\begin{equation}
\label{eq:ineq}
\emph{tr}\big((B_K)_p^2\big)\geq nH_K(p)^2, \quad\text{for any } p\in\Sg.
\end{equation}
Moreover, if $\Sg$ is connected and equality holds for any $p\in\Sg$, then either $H_K=0$ and $\Sg$ is contained in a hyperplane, or $H_K\neq 0$ and $\Sg\subseteq p_0+\la\,(\ptl K)$ for some $p_0\in\rrn$ and $\la>0$.
\end{proposition}

\subsection{Anisotropic area and perimeter}
\label{subsec:perimeter}
\noindent

Along this subsection we represent by $K$ an arbitrary convex body in $\rrn$. Consider a two-sided smooth hypersurface $\Sg\subset\rrn$. Once we have fixed a smooth unit normal vector field $N$ on $\Sg$, the \emph{anisotropic area} is the positive integral
\begin{equation}
\label{eq:area}
\ar_K(\Sg):=\int_\Sg h_K(N)\,dA,
\end{equation}
which is computed with respect to the $n$-dimensional Hausdorff measure. We remark that, when $K$ is not centrally symmetric about $0$, the value of $\ar_K(\Sg)$ may depend on the chosen normal $N$.

On the other hand, the anisotropic boundary area of sets can be introduced by adapting the approach of Caccioppoli and De Giorgi in the isotropic case, see Amar and Bellettini~\cite[Def.~3.1, Prop.~3.2]{amar-bellettini}. For an open set $\Om$ and a Borel set $E$ in $\rrn$, the \emph{anisotropic perimeter of $E$ in $\Om$} is
\begin{equation}
\label{eq:anper}
\pp_K(E,\Om):=\sup\left\{\int_E\divv X\,dV\,;\,X\in C_0^\infty(\Om,K)\right\},
\end{equation}
where $\divv$ is the Euclidean divergence, $V$ denotes the Euclidean volume and $C_0^\infty(\Om,K)$ is the space of smooth vector fields with compact support in $\Om$ and image within $K$. Clearly $\pp_K(E,\Om)$ does not change when $E$ is modified by sets of volume zero. It follows as in \cite[Prop.~3.1]{giusti} that we can always take $E$ so that $0<V(E\cap B)<V(B)$ for any open ball $B\subset\rrn$ centered at $\ptl E\cap\Om$.

If $E$ is a set of locally finite perimeter in $\Om$, then the divergence theorem in \cite[Eq.~(3.47)]{afp} yields
\begin{equation*}
\pp_K(E,\Om)=\sup\left\{\int_{\ptl^*E}\escpr{X,N_E}\,dA\,;\,X\in C_0^\infty(\Om,K)\right\},
\end{equation*}
where $\ptl^*E$ and $N_E$ stand for the reduced boundary and the measure-theoretic outer unit normal of $E$ in $\Om$, respectively. Thus, we infer from \eqref{eq:hK} that $\pp_K(E,\Om)\leq\int_{\ptl^*E}h_K(N_E)\,dA$, where the last integral is the $h_K$-surface energy introduced in \cite[Eq.~(20.2)]{maggi}. Indeed, we have
\begin{equation}
\label{eq:perfinite}
\pp_K(E,\Om)=\int_{\ptl^*E}h_K(N_E)\,dA,
\end{equation}
see Amar and Bellettini~\cite[p.~110]{amar-bellettini}. Hence, when $E$ is open and $\ptl E\cap\Om$ coincides, up to a closed subset of vanishing area, with a smooth hypersurface $\Sg$, we deduce
\begin{equation}
\label{eq:persmooth}
\pp_K(E,\Om)=\int_{\Sg}h_K(N)\,dA=\ar_K(\Sg),
\end{equation}
where $N$ is the outer unit normal on $\Sg$. We note that, for a smooth strictly convex body $K$, equality \eqref{eq:persmooth} comes as in the proof of \cite[Eq.~(1.1)]{giusti} by replacing the unit normal $N$ with the anisotropic normal $N_K$ defined in \eqref{eq:normal}.

When $K$ is the round unit ball about $0$ we get that $\ar_K(\Sg)$ and $\pp_K(E,\Om)$ coincide with the Euclidean area $\ar(\Sg)$ and perimeter $\pp(E,\Om)$, respectively. For an arbitrary $K$, the continuity of $h_K$ on $\sph^n$ implies the existence of constants $\alpha,\beta>0$ depending only on $K$, such that
\begin{equation}
\label{eq:equiv}
\alpha\,\pp(E,\Om)\leq\pp_K(E,\Om)\leq\beta\,\pp(E,\Om),
\end{equation}
for any open set $\Om$ and any set $E$ of locally finite perimeter in $\Om$.

The next result gathers the main properties of the anisotropic perimeter needed in this work.

\begin{proposition}
\label{prop:perprop}
Let $K$ be a convex body and $\Om$ an open set in $\rrn$.
\begin{itemize}
\item[(i)] $($Lower semicontinuity$)$. If $\{E_i\}_{i\in\nn}$ is a sequence of Borel sets in $\Om$ and $\{E_i\}\to E$ in $L^1_{loc}(\Om)$, then $\pp_K(E,\Om)\leq\liminf_{i\to\infty}\pp_K(E_i,\Om)$.
\item[(ii)] $($Compactness$)$. Suppose that $\Om$ is bounded and has Lipschitz boundary. If $\{E_i\}_{i\in\nn}$ is a sequence of Borel sets in $\Om$ such that $\{\pp_K(E_i,\Om)\}_{i\in\nn}$ is bounded, then there is a Borel set $E$ in $\Om$ and a subsequence $\{E_{i(k)}\}_{k\in\nn}$ such that $\{E_{i(k)}\}\to E$ in $L^1(\Om)$.
\item[(iii)] $($Behaviour with respect to translations and dilations$)$. If $E$ is a Borel set in $\Om$, then
\begin{align*}
\pp_K(p_0+E,p_0+\Om)&=\pp_K(E,\Om),
\\
\pp_K(\lambda\,E,\lambda\,\Om)&=\la^n\,\pp_K(E,\Om),
\end{align*}
where $p_0\in\rrn$ and $\lambda>0$.
\item[(iv)] Given a set $E$ of locally finite perimeter in $\Om$ and a point $p\in\rrn$, the inequality
\[
\pp_K(E\setminus B(p,r),\Om)\leq\pp_K(E,\Om\setminus B(p,r))+\ar_K(E\cap\ptl B(p,r))
\]
holds for almost every $r>0$, where $B(p,r)$ is the corresponding open ball, and the last area is computed with respect to the unit normal pointing inside $B(p,r)$.
\end{itemize}
\end{proposition}

\begin{proof}
Statement (i) can be derived as in the isotropic case~\cite[Thm.~1.9]{giusti}. Statement (ii) is a consequence of the first estimate in \eqref{eq:equiv} and the analogous result for the isotropic perimeter~\cite[Thm.~1.19]{giusti}. The two equalities in (iii) come directly from \eqref{eq:anper} and the change of variables formula. The proof of (iv) goes as follows. For a given $X\in C^\infty_0(\Om,\rrn)$, and slightly adapting the approximating sequence $g_\eps(s)$ employed in the proof of \cite[Lem.~1 in Sect.~5.7.1]{evans-gariepy}, we get
\[
\int_{E\setminus B(p,r)}\divv X\,dV=\int_{\ptl^*E\cap(\Om\setminus B(p,r))}\escpr{X,N_E}\,dA+\int_{E\cap\ptl B(p,r)}\escpr{X,\mathcal{N}}\,dA,
\]
where $\mathcal{N}$ is the unit normal on $\ptl B(p,r)$ pointing inside $B(p,r)$. This equality holds for almost every $r>0$. When $X(\Om)\subseteq K$ and having in mind \eqref{eq:hK}, \eqref{eq:perfinite} and \eqref{eq:area}, we deduce
\[
\int_{E\setminus B(p,r)}\divv X\,dV\leq\pp_K(E,\Om\setminus B(p,r))+\ar_K(E\cap\ptl B(p,r)),
\]
which entails the desired inequality.
\end{proof}

\section{The anisotropic partitioning problem}
\label{sec:partprob}
\setcounter{equation}{0}

Here we start discussing the problem of minimizing the anisotropic perimeter among sets with fixed volume inside a given domain. We will review existence, regularity and approximation properties for the solutions, as that as basic facts about the associated isoperimetric profile. These will be employed to prove our main results in Section~\ref{sec:main}. We will denote by $K\subset\rrn$ an arbitrary convex body as defined in Section~\ref{subsec:convexity}. Further hypotheses on $K$ will be specified when needed.

Let $\Om\subset\rrn$ be a bounded open set. The \emph{anisotropic isoperimetric profile} of $\Om$ is the function
\begin{equation}
\label{eq:profile}
I_{\Om,K}(v):=\inf\{\pp_K(E,\Om)\,;\, E\subset\Om,\, V(E)=v\}, \quad v\in[0,V(\Om)],
\end{equation}
where we consider Borel sets with $\pp_K(E,\Om)<\infty$. It is clear that $I_{\Om,K}$ vanishes at the extreme points. From \eqref{eq:equiv} we have
\begin{equation}
\label{eq:equivprofile}
\alpha\,I_{\Om}\leq I_{\Om,K}\leq\beta\,I_{\Om},
\end{equation}
where $I_{\Om}$ stands for the isoperimetric profile when $K$ is the unit ball about the origin. By intersecting $\Om$ with round balls about a point of $\Om$, it follows that $I_{\Om,K}$ is bounded. We also have that $I_{\Om,K}>0$ on $(0,V(\Om))$ because $I_\Om>0$ on $(0,V(\Om))$, see for instance \cite[Thm.~9.3, Lem.~3.4]{ritore-book}. Note that the definition of $I_{\Om,K}(v)$ in \eqref{eq:profile} is still valid when $\Om$ is unbounded and $v\in[0,V(\Om))$.

An \emph{anisotropic isoperimetric region} or \emph{anisotropic minimizer} of volume $v\in(0,V(\Om))$ is a set $E\subset\Om$ with $V(\Om)=v$ and $\pp_K(E,\Om)=I_{\Om,K}(v)$. This is equivalent to that $\pp_K(E,\Om)\leq\pp_K(E',\Om)$, for any other set $E'\subset\Om$ with $V(E')=v$. The existence of minimizers when $\Om$ is bounded is guaranteed by the \emph{direct method}. Indeed, the compactness property and the lower semicontinuity of perimeter in Proposition~\ref{prop:perprop} allow us to reason as in \cite[Prop.~12.30, Re.~20.5]{maggi} to deduce this:

\begin{proposition}
\label{prop:existence}
If $\Om\subset\rrn$ is a bounded open set with Lipschitz boundary, then there are anisotropic isoperimetric regions in $\Om$ of any given volume. 
\end{proposition}

This existence result together with properties of the perimeter leads to the next consequence:

\begin{proposition}
\label{prop:conti}
For a bounded open set $\Om\subset\rrn$ with Lipschitz boundary the profile function $I_{\Om,K}:[0,V(\Om)]\to\rr$ is continuous.
\end{proposition}

\begin{proof}
The continuity at the extremes comes from \eqref{eq:profile} and the second estimate in \eqref{eq:equiv} by taking small round balls inside $\Om$ (or their complementaries). The continuity $I_{\Om,K}$ in $(0,V(\Om))$ is obtained from Propositions~\ref{prop:existence} and \ref{prop:perprop} by following \cite[Re.~3.6]{ritore-book}. The argument also relies on the positivity of $I_{\Om,K}$ and the existence of local deformations for sets of positive perimeter, see \cite[Thm.~1.50]{ritore-book}.
\end{proof}

\begin{remark}
\label{re:centrally}
Suppose that the convex body $K$ is centrally symmetric about $0$. In this case $I_{\Om,K}$ is symmetric about $V(\Om)/2$ since $V(\Om\setminus E)=V(\Om)-V(E)$ and $\pp_K(\Om\setminus E,\Om)=\pp_K(E,\Om)$ for any set $E\subset\Om$. This comes directly from \eqref{eq:anper} since, for any $X\in C^\infty_0(\Om,K)$, we have $-X\in C^\infty_0(\Om,K)$ and
\[
\int_E\divv X\,dV=\int_\Om\divv X\,dV-\int_{\Om\setminus E}\divv X\,dV=\int_{\Om\setminus E}\divv (-X)\,dV
\]
by the divergence theorem. In particular, a set $E\subset\Om$ is an anisotropic isoperimetric region of volume $v$ if and only if $\Om\setminus E$ is an anisotropic isoperimetric region of volume $V(\Om)-v$.
\end{remark}

Next, for an anisotropic minimizer in a smooth domain we discuss analytical properties of the interior boundary and estimates for the Hausdorff dimension of its singular set. The interior regularity was analyzed by Almgren~\cite{almgren} and by Schoen, Simon and Almgren~\cite{schoen-simon-almgren}, see also \cite{bombieri} and \cite{schoen-simon}. The regularity along the free boundary for solutions to more general anisotropic variational problems was studied by De Philippis and Maggi~\cite{philippis-maggi-1,philippis-maggi-2} after previous work by Hardt~\cite{hardt}, and by Duzaar and Steffen~\cite{duzaar-steffen}. These authors considered parametric elliptic integrands depending on the position and the normal direction. Such integrands must satisfy certain conditions which are fulfilled in our setting when $K$ is a smooth strictly convex body. As a consequence of their results we get the following:

\begin{proposition}
\label{prop:regularity}
Let $K\subset\rrn$ be a smooth strictly convex body. Suppose that $\Om\subset\rrn$ is a smooth open set and $E\subset\Om$ is an anisotropic isoperimetric region. Then, the interior boundary $\overline{\ptl E\cap\Om}$ is the disjoint union of a smooth hypersurface $\Sg$, possibly with boundary $\ptl\Sg=\Sg\cap\ptl\Om$, and a relatively closed set of singularities $\Sg_0$ with $\mathcal{H}^{n-2}(\Sg_0)=0$. 
\end{proposition}

\begin{remarks}
\label{re:sentinel1}
(i). Since the Hausdorff measure $\mathcal{H}^{n-2}$ vanishes on $\Sg_0$, the interior boundary has no singularities when $n\leq 2$. In the isotropic case the regularity theory entails that $\Sg_0$ is empty when $n\leq 6$, whereas $\mathcal{H}^{n-7+\gamma}(\Sg_0)=0$ for any $\gamma>0$ when $n\geq 7$. We cannot expect this to hold for anisotropic functionals, as Morgan showed in \cite{morgan-cone-regular} that the cone in $\rr^4$ over the Clifford torus $\sph^1\times\sph^1$ minimizes the perimeter $\pp_K$ for some smooth strictly convex body $K$.

(ii). As indicated in \cite[Re.~1.6]{philippis-maggi-1}, since an anisotropic minimizer verifies an Euler-Lagrange equation in weak form, we can use Schauder's theory for second order elliptic equations to improve the a priori $C^{1,\alpha}$ regularity of $\Sg$ to the $C^\infty$ regularity stated in the proposition, see \cite[Ch.~27]{maggi} for details.

(iii). The interior regularity holds even if $\Om$ is not smooth. Hence, for any anisotropic minimizer $E$ in $\Om$ we have $\ptl E\cap\Om=\Sg\cup\Sg_0$, where $\Sg$ is a smooth hypersurface and $\mathcal{H}^{n-2}(\Sg_0)=0$.

(iv). By taking into account the previous proposition and that the anisotropic perimeter  does not change by sets of volume zero, we can always assume that an anisotropic minimizer $E$ in $\Om$ is an open set. Moreover, equation \eqref{eq:persmooth} implies that $\pp_K(E,U)=\ar_K(\Sg\cap U)$ for any open set $U\subseteq\Om$   because $E$ has almost smooth interior boundary.
\end{remarks}

In the isotropic case, the $\mathcal{H}^{n-2}$-negligibility of the singular set together with upper density estimates leads to find suitable approximations of the function $1$ which are useful to deduce properties of the minimizers from the second derivative of the area. In the next statement we produce similar sequences in the anisotropic setting. An analogous construction is given in the proof of \cite[Lem.~2.5]{philippis-maggi-2}.

\begin{proposition}
\label{prop:approx}
Let $K\subset\rrn$ be a smooth strictly convex body. Suppose that $\Om\subset\rrn$ is a smooth open set, $E\subset\Om$ is an anisotropic minimizer, and $\Sg$ is the regular part of $\overline{\ptl E\cap\Om}$. Then, there is a sequence $\{\omega_\eps\}_{\eps>0}\subset C^\infty_0(\Sg)$ such that:
\begin{itemize}
\item[(i)] $0\leq\omega_\eps\leq 1$ and $\omega_\eps\neq 0$ for any $\eps>0$,
\item[(ii)] $\{\omega_\eps(p)\}\to 1$ for any $p\in\Sg$ when $\eps\to 0$,
\item[(iii)] $\lim_{\eps\to 0}\int_\Sg|\nabla_\Sg\omega_\eps|^2\,dA=0$, where $\nabla_\Sg$ is the gradient in $\Sg$ for the induced metric.
\end{itemize}
\end{proposition}

\begin{proof}
The functions $\omega_\eps$ can be defined as in the isotropic case by following the idea of Sternberg and Zumbrun~\cite[Lem.~2.4]{sz}. Though they assumed $\Om$ to be bounded (and so, the interior boundary $\overline{\ptl E\cap\Om}$ is compact), their argument is extended to the general case, see the proof of \cite[Thm.~3.6]{rosales-gauss} for the details. The main ingredients are the fact that the singular set $\Sg_0$ of $\overline{\ptl E\cap\Om}$ satisfies $\mathcal{H}^{n-2}(\Sg_0)=0$, and the existence of constants $C_0,R_0>0$ such that 
\begin{equation}
\label{eq:densesti}
\ar(\Sg\cap B(p,r))\leq C_0\,r^n,
\end{equation}
for any open ball $B(p,r)\subset\rr^n$ with $r\leq R_0$. As the first condition holds by Proposition~\ref{prop:regularity}, it suffices to show that the anisotropic area minimality of $E$ entails the upper density estimate in \eqref{eq:densesti}. This kind of estimates are well-known for minimizers of perimeter functionals, see for instance \cite[Thm.~16.14]{maggi} and \cite[Eq.~(2.47)]{philippis-maggi-1}. For the sake of completeness we give here a quick proof.

Choose $p_0\in\rrn$ and $R_0>0$ for which $\overline{B}(p_0,R_0)\subset\Om\setminus\overline{E}$. Fix a point $p\in\rrn$ and consider a ball $B(p,r)$ with $r\leq R_0$ and $E\cap B(p,r)\neq\emptyset$. Since $V(E\cap B(p,r))\leq V(B(p,r))\leq V(B(p_0,R_0))$, there is $r'>0$ such that $V(B(p_0,r'))=V(E\cap B(p,r))$. It is clear that $r'\leq r$. By setting
\[
E':=(E\setminus B(p,r))\cup\overline{B}(p_0,r'),
\] 
we produce a set in $\Om$ with $V(E')=V(E)$. By having in mind that $\pp_K(E,\cdot)$ defines a finite Borel measure in $\Om$, the isoperimetry of $E$, and the inequality in Proposition~\ref{prop:perprop} (iv), we obtain
\begin{align*}
\pp_K(E,\Om\cap B(p,r))&+\pp_K(E,\Om\setminus B(p,r))\leq\pp_K(E,\Om)\leq\pp_K(E',\Om)
\\
&=\pp_K(E\setminus B(p,r),\Om)+\pp_K(B(p_0,r'))
\\
&\leq\pp_K(E,\Om\setminus B(p,r))+\ar_K(E\cap\ptl B(p,r))+\pp_K(B(p_0,r')).
\end{align*}
This comparison holds for almost every $r\leq R_0$. By simplifying the term $\pp_K(E,\Om\setminus B(p,r))$, taking into account  \eqref{eq:persmooth} and \eqref{eq:equiv}, and using that $r'\leq r$, we have
\begin{align*}
\pp_K(E,\Om\cap B(p,r))&\leq\pp_K(B(p,r))+\pp_K(B(p_0,r'))
\\
&\leq\beta\,\big(\pp(B(p,r))+\pp(B(p_0,r'))\big)\leq 2\beta\,\la_n\,r^n,
\end{align*}
where $\la_n:=\pp(B(0,1))$. Hence, for a fixed $p\in\rrn$, we can employ the first inequality in \eqref{eq:equiv} and the fact that $\pp(E,\Om\cap B(p,r))=\ar(\Sg\cap B(p,r))$ to deduce the estimate in \eqref{eq:densesti} for almost every $r>0$. The validity of \eqref{eq:densesti} for every $r>0$ follows by the dominated convergence theorem.
\end{proof}

We finish this section by reviewing the anisotropic partitioning problem in the special case where the ambient domain is a cone. The explicit expression of the isoperimetric profile and the description of certain solutions will be useful for some of our results. 

Let $K\subset\rrn$ be a convex body and $\cc\subseteq\rrn$ an open cone with vertex at the origin (we do not assume additional regularity on $K$ or $\cc$). Since $\cc$ and $\ptl\cc$ are invariant under dilations about $0$ we can reason as in \cite[Prop.~3.1]{cones} with the help of Proposition~\ref{prop:perprop} (iii) to infer the existence of a constant $c(\cc,K)\geq 0$ such that 
\begin{equation}
\label{eq:isocones}
I_{\cc,K}(v)=c(\cc,K)\,v^{n/(n+1)},\quad\text{for any } v\geq 0.
\end{equation} 
In particular $(I_{\cc,K})^{(n+1)/n}$ is a linear function.

When $\cc$ is also convex, the sets $\la K\cap\cc$ provide anisotropic isoperimetric regions in $\cc$ for any $\la>0$, see  Cabr\'e, Ros-Oton and Serra~\cite[Thm.~1.3]{cabre-rosoton-serra}, and Milman and Rotem~\cite[Cor.~1.2]{milman-rotem}, after previous work of Lions and Pacella in the isotropic case~\cite[Thm.~1.1]{lp}. Hence, we get
\[
c(\cc,K)=\frac{\pp_K(K,\cc)}{V(K\cap\cc)^{n/(n+1)}}.
\]
On the other hand, we have the following relation (which holds for any Lipschitz cone, convex or not)
\begin{equation}
\label{eq:relationcone}
\pp_K(K,\cc)=(n+1)\,V(K\cap\cc).
\end{equation}
This comes from the divergence theorem applied to the identity vector field over $K\cap\cc$, see \cite[Eq.~(1.14)]{cabre-rosoton-serra} for the details. From here we deduce that
\begin{equation}
\label{eq:isoprofcones}
I_{\cc,K}(v)=(n+1)\,V(K\cap\cc)^{1/(n+1)}\,v^{n/(n+1)}, \quad\text{for any } v\geq 0.
\end{equation}
As a consequence $I_{\cc_1,K_1}\leq I_{\cc_2,K_2}$ whenever $V(K_1\cap\cc_1)\leq V(K_2\cap\cc_2)$. 

\begin{remark}
\label{re:halfspace}
The aforementioned results apply in an open half-space $\hh$ with $0\in\ptl\hh$. Unlike the isotropic case, the profile $I_{\hh,K}$ may depend on $\hh$ because the volume $V(K\cap\hh)$ could be a non-constant function of $\hh$. When the convex body $K$ is centrally symmetric about $0$ we know that $V(K\cap\hh)=V(K)/2$, and so $I_{\hh,K}$ does not depend on the half-space $\hh$.
\end{remark}

\section{Variational formulas and consequences}
\label{sec:2ndvar}
\setcounter{equation}{0}

We begin this section by recalling the first variation of the anisotropic area and the equilibrium conditions satisfied by the solutions to the anisotropic partitioning problem. Then, we compute the second derivative of the area for certain deformations. From here, and taking advantage of the equilibrium conditions, we will infer a second variation formula that will play a central role in Section~\ref{sec:main}.  

Let $K\subset\rrn$ be a convex body and $\Sg\subset\rrn$ a two-sided smooth hypersurface, possibly with boundary $\ptl\Sg$. We will follow the notation introduced in Section~\ref{subsec:extrinsic} for the anisotropic extrinsic geometry of $\Sg$. In the case $\ptl\Sg=\emptyset$ we will assume that all the terms and integrals over $\ptl\Sg$ vanish. 

Consider a complete smooth vector field $X$ on $\rrn$ with one-parameter group of diffeomorphisms $\{\phi_t\}_{t\in\rr}$. The associated \emph{variation} of $\Sg$ is the family $\{\Sg_t\}_{t\in\rr}$, where $\Sg_t:=\phi_t(\Sg)$ for any $t\in\rr$. Once we have fixed a smooth unit normal $N$ on $\Sg$, we can find a smooth vector field $\ovn$, whose restriction to $\Sg_t$ provides a unit normal $N_t$ with $N_0=N$. If $X$ has compact support on $\Sg$ then there is a compact set $C\subseteq\Sg$ such that $\phi_t(p)=p$ for any $p\in\Sg\setminus C$ and $t\in\rr$. The corresponding \emph{anisotropic area functional} is defined as
\begin{equation}
\label{eq:areafunct}
\ar_K(t):=\ar_K(\phi_t(C))=\int_{\phi_t(C)}h_K(N_t)\,dA_t=\int_C h_K(\ovn\circ\phi_t)\,\text{Jac}\,\phi_t\,dA,
\end{equation}
where $\text{Jac}\,\phi_t$ is the Jacobian of the diffeomorphism $\phi_{t|\Sg}:\Sg\to\Sg_t$. 

The computations of $\ar_K'(0)$ and $\ar_K''(0)$ are key tools to understand the critical points and the second order minima of the anisotropic area. When $K$ is smooth and strictly convex the first variation is well-known, see for instance Clarenz~\cite[Sect.~1]{clarenz-2002}, Koiso and Palmer~\cite[Proof of Prop.~3.1]{koiso-palmer2} or Koiso~\cite[Lem.~9.1]{koiso-2019}. By following \cite[Prop.~3.2]{rosales-cones} we get
\begin{equation}
\label{eq:1stareagen}
\ar_K'(0)=-\int_\Sg nH_K\,u\,dA-\int_{\ptl\Sg}\escpr{X,\nu_K}\,dL,
\end{equation}
where $u:=\escpr{X,N}$ and $L$ stands for the $(n-1)$-dimensional Hausdorff measure. On the other hand, from \cite[\S9]{simon} we know that, if $E\subset\rrn$ is a set of locally finite perimeter with $V(E)<\infty$, and $X$ is a smooth vector field with compact support on $\rrn$, then the first derivative of the \emph{volume functional} $V(t):=V(E_t)$, where $E_t:=\phi_t(E)$ for any $t\in\rr$, is given by
\begin{equation}
\label{eq:1stvolgen}
V'(0)=\int_E\divv X\,dV=\int_{\ptl^*E}\escpr{X,N_E}\,dA,
\end{equation}
where $\ptl^*E$ is the reduced boundary and $N_E$ the measure-theoretic outer unit normal. 

The next result gathers the equilibrium conditions (in strong form) of an anisotropic isoperimetric region $E$ within a smooth open set $\Om$. These come from the formulas \eqref{eq:1stareagen} and \eqref{eq:1stvolgen} because the regular part $\Sg$ of $\overline{\ptl E\cap\Om}$ is critical for the area $\mathcal{A}_K$ (computed with respect to the outer unit normal) under volume-preserving deformations moving compact portions of $\Sg$ while leaving invariant $\ptl\Om$, see Koiso and Palmer~\cite[Prop.~3.1]{koiso-palmer2}, and De Philippis and Maggi~\cite[Eqs.~(1.4), (1.5), Thm.~1.2]{philippis-maggi-1}.

\begin{proposition}
\label{prop:varprop}
Let $K\subset\rrn$ be a smooth strictly convex body. Suppose that $\Om\subset\rrn$ is a smooth open set, $E\subset\Om$ is an anisotropic minimizer, and $\Sg$ is the regular part of $\overline{\ptl E\cap\Om}$. Then, the anisotropic mean curvature $H_K$ of $\Sg$ is constant with respect to the outer unit normal $N$. Moreover, the anisotropic normal $N_K$ and the inner unit normal $\xi$ to $\ptl\Om$ satisfy $\escpr{N_K,\xi}=0$ along $\ptl\Sg$.
\end{proposition}

\begin{remark}
\label{re:sentinel2}
When $\Om$ is not smooth we can still deduce that the regular part of $\ptl E\cap\Om$ has constant anisotropic mean curvature with respect to the outer unit normal.
\end{remark}

The equilibrium conditions have many consequences. For instance, it was shown by the author~\cite[Eq.~(4.1), Lem.~4.3]{rosales-cones} that the equality $\escpr{N_K,\xi}=0$ along $\ptl\Sg$ entails
\begin{equation}
\label{eq:conormalest}
\escpr{\nu,\xi}(p)\neq 0\quad\text{and}\quad\nu_K(p)=\frac{\var_K}{\escpr{\nu,\xi}}(p)\,\xi(p),\quad\text{for any } p\in\ptl\Sg,
\end{equation}
where $\nu$ is the inner conormal and $\nu_K$ the associated anisotropic conormal. Thus, the boundary term in \eqref{eq:1stareagen} vanishes when $X$ is tangent to $\ptl\Om$. This allows to simplify the first variation formulas for some special deformations that will be mostly employed in this work.

\begin{corollary}
\label{cor:1stvarest}
In the conditions of Proposition~\ref{prop:varprop}, if $X$ is a smooth vector field with compact support on $\rrn$ such that $X$ is tangent to $\ptl\Om$ and $X_{|\Sg}=\omega N_K$ for some $\omega\in C^\infty_0(\Sg)$, then we have
\[
\ar_K'(0)=-nH_K\int_\Sg\omega\,\var_K\,dA \quad\text{and}\quad V'(0)=\int_\Sg\omega\,\var_K\,dA.
\]
In particular $(\ar_K+nH_KV)'(0)=0$.
\end{corollary}

Now, we establish an expression for $\ar_K''(0)$ when an arbitrary hypersurface $\Sg$ is moved under the flow of a vector field $X$ as in the previous corollary. For the proof we will follow the scheme in \cite[Prop.~3.3]{rosales-cones}, where we discussed the special case $\omega=1$ for a compact hypersurface $\Sg$. 

\begin{proposition}
\label{prop:2ndarea}
Let $K\subset\rrn$ be a smooth strictly convex body and $\Sg\subset\rrn$ a two-sided smooth hypersurface, possibly with boundary $\ptl\Sg$, and unit normal $N$. Take a complete smooth vector field $X$ on $\rrn$ such that $X_{|\Sg}=\omega N_K$ for some $\omega\in C^\infty_0(\Sg)$. Then, we have
\begin{align*}
\ar_K''(0)&=\int_\Sg\escpr{\nabla^K_\Sg\omega,\nabla_\Sg\omega}\,\var^2_K\,dA+\int_\Sg\left(n^2H_K^2-\emph{tr}(B_K^2)\right)\omega^2\var_K\,dA
\\
&-\int_\Sg n H_K\,v\,dA-\int_{\ptl\Sg}\escpr{Z,\nu_K}\,dL,
\end{align*}
where $v:=\escpr{Z,N}$ and $Z:=D_XX$.
\end{proposition}

\begin{proof}
For any $w\in\rrn$ the notations $w^\top$ and $w^\bot$ stand for the projections of $w$ onto $T\Sg$ and $(T\Sg)^\bot$, respectively. For a fixed $p\in\Sg$ we consider the functions
\begin{equation}
\label{eq:hyj}
h_p(t):=h_K(N_t\circ\phi_t)(p), \quad j_p(t):=(\text{Jac}\,\phi_t)(p), \quad\text{for any }t\in\rr.
\end{equation}
It is clear from \eqref{eq:areafunct} that $\ar_K(t)=\int_C h_p(t)\,j_p(t)\,dA$. By differentiating under the integral sign twice and taking into account that $\phi_0=\text{Id}$, we infer
\begin{equation}
\label{eq:2ndaK}
\ar_K''(0)=\int_C\big(h_p''(0)+2h_p'(0)\,j_p'(0)+\var_K(p)\,j_p''(0)\big)\,dA.
\end{equation}
Let us compute all the derivatives in the integrand above.

For the calculus of $j_p'(0)$ and $j_p''(0)$ we refer to Simon~\cite[\S9]{simon}. On the one hand, it is known that $j_p'(0)=(\divv_\Sg X)(p)$. So, from \eqref{eq:mc} and the fact that $X_{|\Sg}=\omega N_K$ we obtain
\begin{equation}
\label{eq:jprima}
j_p'(0)=(\divv_\Sg X)(p)=-nH_K(p)\,\omega(p)+\escpr{N_K,\nabla_\Sg\omega}(p).
\end{equation}
On the other hand
\begin{equation}
\label{eq:j2preprima}
j_p''(0)=(\divv_\Sg Z)(p)+(\divv_\Sg X)^2(p)+\sum_{i=1}^n|(D_{e_i}X)^\bot|^2-\sum_{i,j=1}^n\escpr{D_{e_i}X,e_j}\,\escpr{D_{e_j}X,e_i},
\end{equation}
where $\{e_1,\ldots,e_n\}$ is an orthonormal basis of $T_p\Sg$. From \eqref{eq:shape} and equality $X_{|\Sg}=\omega N_K$ we get
\begin{equation}
\label{eq:pito}
D_{e_i}X=\escpr{(\nabla_\Sg\omega)(p),e_i}\,N_K(p)-\omega(p)\,(B_K)_p(e_i),
\end{equation}
and so $(D_{e_i}X)^\bot=\escpr{(\nabla_\Sg\omega)(p),e_i}\,\var_K(p)\,N(p)$ by \eqref{eq:nkn}. As a consequence
\[
\sum_{i=1}^n|(D_{e_i}X)^\bot|^2=\var_K(p)^2\,|\nabla_\Sg\omega|^2(p).
\]
The expression for $D_{e_i}X$ in \eqref{eq:pito} also yields
\begin{align*}
\sum_{i,j=1}^n\escpr{D_{e_i}X,e_j}\,\escpr{D_{e_j}X,e_i}&=\sum_{i,j=1}^n\escpr{N_K(p),e_i}\,\escpr{N_K(p),e_j}\,\escpr{(\nabla_\Sg\omega)(p),e_i}\,\escpr{(\nabla_\Sg\omega)(p),e_j}
\\
&+\omega(p)^2\,\sum_{i,j=1}^n\escpr{(B_K)_p(e_i),e_j}\,\escpr{(B_K)_p(e_j),e_i}
\\
&-2\,\omega(p)\,\sum_{i,j=1}^n\escpr{(\nabla_\Sg\omega)(p),e_j}\,\escpr{N_K(p),e_i}\,\escpr{(B_K)_p(e_i),e_j}
\\
&=\escpr{N_K,\nabla_\Sg\omega}(p)^2+\omega(p)^2\,\text{tr}((B_K)_p^2)
\\
&-2\,\omega(p)\,\sum_{i=1}^n\escpr{N_K(p),e_i}\,\escpr{(\nabla_\Sg\omega)(p),(B_K)_p(e_i)}.
\end{align*}
By substituting \eqref{eq:jprima} and the previous calculus into \eqref{eq:j2preprima}, and simplifying, we arrive at
\begin{equation}
\label{eq:j2prima}
\begin{split}
j_p''(0)&=\big(\!\divv_\Sg Z+n^2H_K^2\,\omega^2+\var_K^2\,|\nabla_\Sg\omega|^2-\omega^2\,\text{tr}(B_K^2)-2nH_K\,\omega\,\escpr{N_K,\nabla_\Sg\omega}\big)(p)
\\
&+2\,\omega(p)\,\sum_{i=1}^n\escpr{N_K(p),e_i}\,\escpr{(\nabla_\Sg\omega)(p),(B_K)_p(e_i)}.
\end{split}
\end{equation}

Next, we compute $h_p'(0)$ and $h_p''(0)$. From \eqref{eq:hyj} and \eqref{eq:gradhK} we deduce  
\begin{equation}
\label{eq:catrina}
h_p'(t)=\escpr{\pi_K((N_t\circ\phi_t)(p)),(D_{X}\ovn)(\phi_t(p))}.
\end{equation}
This implies that
\begin{equation}
\label{eq:hprima0}
h_p'(0)=\escpr{N_K(p),D_{X(p)}\overline{N}}.
\end{equation}
The calculus of $D_{X(p)}\ovn$ is found in \cite[Lem.~4.1(1)]{ros-souam} and involves the function $u:=\escpr{X,N}$. In our case $u=\escpr{\omega N_K,N}=\omega\,\var_K$ by \eqref{eq:nkn}. So, we get
\begin{equation}
\label{eq:tente}
\begin{split}
D_{X(p)}\ovn&=\frac{d}{dt}\bigg|_{t=0}(N_t\circ\phi_t)(p)=-(\nabla_\Sg u)(p)-B_p(X^\top(p))
\\
&=-\omega(p)\,(\nabla_\Sg\var_K)(p)-\var_K(p)\,(\nabla_\Sg\omega)(p)-\omega(p)\,B_p(N_K^\top(p))
\\
&=-\var_K(p)\,(\nabla_\Sg\omega)(p),
\end{split}
\end{equation}
where in the last equality we have employed the identity
\begin{equation}
\label{eq:gradvarK}
\nabla_\Sg\var_K=-B(N_K^\top),
\end{equation}
which was proved in \cite[Lem.~3.1]{rosales-cones}. By substituting \eqref{eq:tente} into \eqref{eq:hprima0} we obtain
\begin{equation}
\label{eq:hprima}
h_p'(0)=-\var_K(p)\,\escpr{N_K,\nabla_\Sg\omega}(p).
\end{equation}
On the other hand, by differentiating at $t=0$ in \eqref{eq:catrina} and taking into account \eqref{eq:endoq}, we arrive at
\[
h_p''(0)=\escpr{Q_p(D_{X(p)}\ovn),D_{X(p)}\ovn}+\escpr{N_K(p),D_{X(p)}D_X\ovn}.
\]
From \eqref{eq:tente} and the formula in Lemma~\ref{lem:potato} below, it follows that
\begin{equation}
\label{eq:h2prima}
\begin{split}
h_p''(0)&=\big(\var_K^2\,\escpr{\nabla_\Sg^K\omega,\nabla_\Sg\omega}+2\,\var_K\,\escpr{N_K,\nabla_\Sg\omega}^2-\escpr{N_K,\nabla_\Sg v}-\var_K^3\,|\nabla_\Sg\omega|^2\big)(p)
\\
&+\escpr{\nabla_\Sg\var_K,Z}(p)-2\,\omega(p)\,\var_K(p)\,\sum_{i=1}^n\escpr{(\nabla_\Sg\omega)(p),(B_K)_p(e_i)}\,\escpr{N_K(p),e_i},
\end{split}
\end{equation}
where we have used~\eqref{eq:gradvarK}, the definition in \eqref{eq:angrad} and the symmetry of the shape operator.

Now, by having in mind \eqref{eq:h2prima}, \eqref{eq:hprima}, \eqref{eq:jprima}, \eqref{eq:j2prima}, and simplifying, we conclude that the integrand in \eqref{eq:2ndaK} is the evaluation at $p$ of the function
\[
\escpr{\nabla_\Sg^K\omega,\nabla_\Sg\omega}\,\var_K^2-\escpr{N_K,\nabla_\Sg v}+\divv_\Sg(\var_K\,Z)+\big(n^2H^2_K-\text{tr}(B_K^2)\big)\,\omega^2\var_K.
\]
From here the proof finishes by applying the formula
\[
\int_\Sg\divv_\Sg(\var_K\,Z)\,dA=-\int_\Sg nH_K\,v\,dA+\int_\Sg\escpr{N_K,\nabla_\Sg v}\,dA-\int_{\ptl\Sg}\escpr{Z,\nu_K}\,dL,
\]
which comes from \cite[Lem.~3.1]{rosales-cones} by taking $U=Z$.
\end{proof}

\begin{lemma}
\label{lem:potato}
In the conditions of Proposition~\ref{prop:2ndarea}, for any $p\in\Sg$, we have
\begin{align}
D_{X(p)}D_X\ovn&=2\,\var_K(p)\,\escpr{N_K,\nabla_\Sg\omega}(p)\,(\nabla_\Sg\omega)(p)-(\nabla_\Sg v)(p)-B_p(Z^\top(p))
\\
&-2\,\omega(p)\,\var_K(p)\,\sum_{i=1}^n\escpr{(\nabla_\Sg\omega)(p),(B_K)_p(e_i)}\,e_i-\var_K(p)^2\,|\nabla_\Sg\omega|^2(p)\,N(p).
\end{align}
\end{lemma}

\begin{proof}
Let $\{e_1,\ldots,e_n\}$ be an orthonormal basis of $T_p\Sg$. By using the flow $\{\phi_t\}_{t\in\rr}$ of $X$ we can construct, for any $i=1,\ldots,n$, a smooth vector field $E_i$ around $p$ which is tangent on any $\Sg_t=\phi_t(\Sg)$ while satisfying $E_i(p)=e_i$ and $[X,E_i]=0$ (here $[\cdot\,,\cdot]$ is the Lie bracket of vector fields). Clearly
\begin{equation}
\label{eq:gen0}
D_{X(p)}D_X\ovn=\sum_{i=1}^n\escpr{D_{X(p)}D_X\ovn,e_i}\,e_i
+\escpr{D_{X(p)}D_X\ovn,N(p)}\,N(p).
\end{equation}
Now we compute the different terms in the previous equation.

Since $|\ovn|^2=1$ we get $\escpr{D_X\ovn,\ovn}=0$. By differentiating and applying equation \eqref{eq:tente}, we deduce
\begin{equation}
\label{eq:gen1}
\escpr{D_{X(p)}D_X\ovn,N(p)}=-|D_{X(p)}\ovn|^2=-\var_K(p)^2\,|\nabla_\Sg\omega|^2(p).
\end{equation}
Next, we differentiate twice with respect to $X$ in equality $\escpr{\ovn,E_i}=0$. We obtain
\begin{equation}
\label{eq:gen01}
\escpr{D_{X(p)}D_X\ovn,e_i}=-2\,\escpr{D_{X(p)}\ovn,D_{e_i}X}-\escpr{N(p),D_{X(p)}D_{E_i}X},
\end{equation}
were we have used that $[X,E_i]=0$. From \eqref{eq:tente} and \eqref{eq:pito} we get
\[
\escpr{D_{X(p)}\ovn,D_{e_i}X}=-\var_K(p)\,\escpr{(\nabla_\Sg\omega)(p),e_i}\,\escpr{N_K,\nabla_\Sg\omega}(p)+\omega(p)\,\var_K(p)\,\escpr{(\nabla_\Sg\omega)(p),(B_K)_p(e_i)}.
\]
On the other hand, as the Riemann curvature tensor vanishes in $\rrn$, we infer
\[
0=D_XD_{E_i}X-D_{E_i}D_XX-D_{[X,E_i]}X=D_XD_{E_i}X-D_{E_i}Z
\]
because $D_XX=Z$. This shows that $D_{X(p)}D_{E_i}X=D_{e_i}Z$. As a consequence
\begin{align}
\escpr{N(p),D_{X(p)}D_{E_i}X}&=\escpr{N(p),D_{e_i}Z}=e_i\big(\escpr{Z,N}\big)-\escpr{Z(p),D_{e_i}N}
\\
&=\escpr{(\nabla_\Sg v)(p),e_i}+\escpr{Z^\top(p),B_p(e_i)}
\\
&=\escpr{(\nabla_\Sg v)(p),e_i}+\escpr{B_p(Z^\top(p)),e_i}.
\end{align}
Coming back to \eqref{eq:gen01} we deduce
\begin{equation}
\label{eq:gen2}
\begin{split}
\escpr{D_{X(p)}D_X\ovn,e_i}&=2\,\var_K(p)\,\escpr{(\nabla_\Sg\omega)(p),e_i}\,\escpr{N_K(p),(\nabla_\Sg\omega)(p)}
\\
&-2\,\omega(p)\,\var_K(p)\,\escpr{(\nabla_\Sg\omega)(p),(B_K)_p(e_i)}-\escpr{(\nabla_\Sg v)(p),e_i}-\escpr{B_p(Z^\top(p)),e_i}.
\end{split}
\end{equation}
The proof finishes by substituting \eqref{eq:gen1} and \eqref{eq:gen2} into \eqref{eq:gen0}.
\end{proof}

In Corollary~\ref{cor:1stvarest} we saw that the equilibrium conditions for solutions to the anisotropic partitioning problem imply the equality $(\ar_K+nH_KV)'(0)=0$ for certain deformations. In the next result we obtain a useful expression for $(\ar_K+nH_KV)''(0)$. We remark that the formula remains valid for sets of locally finite perimeter satisfying the regularity properties in Proposition~\ref{prop:regularity} (we only need $\ar(\Sg_0)=0$) and the equilibrium conditions in Proposition~\ref{prop:varprop}. Previous second variational formulas for stationary hypersurfaces in anisotropic capillarity problems were established by Koiso and Palmer~\cite[Prop.~3.3]{koiso-palmer2}, \cite[Prop.~3.3]{koiso-palmer4}, De Philippis and Maggi~\cite[Lem.~A.5]{philippis-maggi-2}, Barbosa and Carvalho Silva~\cite[Prop.~3]{barbosa-silva}, and Guo and Xia~\cite[Prop.~3.5]{guo-xia}. 

\begin{proposition}
\label{prop:2ndareavol}
Let $K\subset\rrn$ be a smooth strictly convex body. Suppose that $\Om\subset\rrn$ is a smooth open set, $E\subset\Om$ is an anisotropic minimizer, and $\Sg$ is the regular part of $\overline{\ptl E\cap\Om}$. If $X$ is a smooth vector field with compact support on $\rrn$ such that $X$ is tangent to $\ptl\Om$ and $X_{|\Sg}=\omega N_K$ for some $\omega\in C^\infty_0(\Sg)$, then we have
\[
(\ar_K+nH_KV)''(0)=\ind_K(\omega).
\]
Here $\ind_K$ is the quadratic form on $C^\infty_0(\Sg)$ defined by
\begin{equation}
\label{eq:indexform}
\ind_K(\omega):=\int_\Sg\left\{\escpr{\nabla^K_\Sg\omega,\nabla_\Sg\omega}\,\var^2_K-\emph{tr}(B_K^2)\,\omega^2\,\var_K\right\}dA-\int_{\ptl\Sg}\frac{\emph{II}(N_K,N_K)}{\escpr{\nu,\xi}}\,\omega^2\,\var_K\,dL,
\end{equation}
where $\emph{II}$ is the second fundamental form of $\ptl\Om$ with respect to $\xi$, and $\nu$ is the inner conormal of $\ptl\Sg$.
\end{proposition}

\begin{remark}
The equality $\escpr{N_K,\xi}=0$ in Proposition~\ref{prop:varprop} means that $N_K$ is tangent to $\ptl\Om$ along $\ptl\Sg$. This guarantees existence of vector fields in the hypotheses of the statement and that $\text{II}(N_K,N_K)$ is well defined. We also note that $\escpr{\nu,\xi}$ never vanishes along $\ptl\Sg$ by equation~\eqref{eq:conormalest}, so that $\Sg$ intersects $\ptl\Om$ transversally along $\ptl\Sg$. 
\end{remark}

\begin{proof}[Proof of Proposition~\ref{prop:2ndareavol}]
First, we compute the second derivative of the volume functional $V(t)$. From Sternberg and Zumbrun~\cite[Eq.~(2.17)]{sz2} we have
\begin{equation}
\label{eq:2ndvol}
V''(0)=\int_E\divv\big((\divv X)\,X\big)\,dV=\int_\Sg(\divv X)\,\escpr{X,N}\,dA=\int_\Sg(\divv X)\,\omega\,\var_K\,dA,
\end{equation}
where we used the divergence theorem for finite perimeter sets, the fact that $X$ is tangent to $\ptl\Om$, the negligibility of $\Sg_0$ for the area measure, and that $X_{|\Sg}=\omega N_K$. By taking into account \eqref{eq:mc} we get
\begin{align*}
(\divv X)\,\omega\,\var_K&=\big(\!\divv_\Sg(\omega N_K)+\escpr{D_NX,N}\big)\,\omega\,\var_K
\\
&=-nH_K\,\omega^2\,\var_K+\escpr{\nabla_\Sg\omega,N_K}\,\omega\,\var_K+\escpr{D_NX,N}\,\omega\,\var_K.
\end{align*}
On the other hand
\begin{align*}
v:=\escpr{Z,N}=\escpr{D_XX,N}=\escpr{D_{X^\top}X+D_{X^\bot}X,N}=\omega\,\escpr{D_{N_K^\top}X,N}+\omega\,\var_K\,\escpr{D_NX,N},
\end{align*}
and so
\[
(\divv X)\,\omega\,\var_K=-nH_K\,\omega^2\,\var_K+\escpr{\nabla_\Sg\omega,N_K}\,\omega\,\var_K+v-\omega\,\escpr{D_{N_K^\top}X,N}.
\]
Observe also that
\begin{align*}
\escpr{D_{N_K^\top}X,N}&=N_K^\top(\escpr{X,N})-\escpr{X,D_{N_K^\top}N}  =\escpr{\nabla_\Sg(\omega\,\var_K),N_K}-\escpr{X,D_{N_K^\top}N}
\\
&=\omega\,\escpr{\nabla_\Sg\var_K,N_K}+\var_K\,\escpr{\nabla_\Sg\omega,N_K}+\escpr{\omega\,N_K,B(N_K^\top)}=\var_K\,\escpr{\nabla_\Sg\omega,N_K},
\end{align*}
where in the last equality we employed \eqref{eq:gradvarK}. Together with the previous equations this shows that
\[
V''(0)=\int_\Sg(-nH_K\,\omega^2\var_K+v)\,dA.
\]
From this expression, the formula in Proposition~\ref{prop:2ndarea}, and the fact that $H_K$ is constant, we deduce
\[
(\ar_K+nH_KV)''(0)=\int_\Sg\left\{\escpr{\nabla^K_\Sg\omega,\nabla_\Sg\omega}\,\var^2_K-\text{tr}(B_K^2)\,\omega^2\var_K\right\}dA-\int_{\ptl\Sg}\escpr{Z,\nu_K}\,dL.
\]
As to the boundary term, from \eqref{eq:conormalest} we obtain
\[
\escpr{Z,\nu_K}=\frac{\var_K}{\escpr{\nu,\xi}}\,\escpr{Z,\xi}\quad\text{along }\ptl\Sg.
\]
Observe also that $\escpr{X,\xi}=0$ along $\ptl\Om$ because $X$ is tangent to $\ptl\Om$. By differentiating with respect to $X$, it follows that
\[
0=\escpr{D_XX,\xi}+\escpr{X,D_X\xi}=\escpr{Z,\xi}-\text{II}(X,X),
\]
and so $\escpr{Z,\xi}=\text{II}(X,X)=\text{II}(N_K,N_K)\,\omega^2$ along $\ptl\Sg$. All this leads to
\[
\escpr{Z,\nu_K}=\frac{\text{II}(N_K,N_K)}{\escpr{\nu,\xi}}\,\omega^2\,\var_K,
\]
which completes the proof.
\end{proof}

Finally, we employ the previous proposition to compute the first and second derivatives of a certain power of the relative anisotropic profile associated to a local perturbation of a minimizer.

\begin{corollary}
\label{cor:2ndrelprof}
In the conditions of Proposition~\ref{prop:2ndareavol}, if we also suppose that $V'(0)\neq 0$, and we consider the function $f(v):=\ar_K(V^{-1}(v))^{(n+1)/n}$ defined locally around $v_0:=V(E)$, then
\begin{align*}
f'(v_0)&=-(n+1)\,H_K\,\ar_K(\Sg)^{1/n},
\\
f''(v_0)&=\frac{n+1}{n}\,\ar_K(\Sg)^{1/n}\left\{\frac{nH_K^2}{\ar_K(\Sg)}+\frac{\ind_K(\omega)}{\big(\int_\Sg\omega\,\var_K\,dA\big)^2}\right\}.
\end{align*}
\end{corollary}

\begin{proof}
Note that $f(v)$ is well defined around $v_0$ because $V$ is a local diffeomorphism at $v_0$. A straightforward calculus shows that
\[
f'(v)=\frac{n+1}{n}\,\ar_K(V^{-1}(v))^{1/n}\,\frac{\ar_K'(V^{-1}(v))}{V'(V^{-1}(v))}.
\]
By differentiating at $v=v_0$ we arrive at
\begin{align}
f''(v_0)=\frac{n+1}{n}\,\ar_K(\Sg)^{1/n}\,\left\{\frac{1}{n\,\ar_K(\Sg)}\,\frac{\ar_K'(0)^2}{V'(0)^2}+\frac{\ar_K''(0)-\frac{\ar_K'(0)}{V'(0)}\,V''(0)}{V'(0)^2}\right\}.
\end{align}
The proof finishes by using the equalities in Corollary~\ref{cor:1stvarest} and that $(\ar_K+nH_KV)''(0)=\ind_K(\omega)$.
\end{proof}

\section{Main results}
\label{sec:main}

In this section we study the anisotropic isoperimetric problem inside convex domains. We will establish analytical properties of the associated isoperimetric profile, sharp comparisons, and topological features about the minimizers. The obtained statements extend known facts in the isotropic setting, see the Introduction for a detailed exposition with references. The proofs are mostly inspired in the ones for the isotropic case by employing the tools developed in Sections~\ref{sec:partprob} and \ref{sec:2ndvar}. 

Our first result concerns the concavity of the anisotropic isoperimetric profile. We need a second order derivative criterion to ensure that a function with low regularity is concave.

\begin{lemma}
\label{lem:weak2nd}
Let $f:I\to\rr$ be a lower semicontinuous function over an open interval $I\subseteq\rr$. If the upper second derivative defined as
\[
\overline{D^2}f(x):=\limsup_{t\to 0^+}\frac{f(x+t)+f(x-t)-2f(x)}{t^2}
\]
satisfies $\overline{D^2}f(x)\leq 0$ for any $x\in I$, then $f$ is concave.
\end{lemma}

\begin{proof}
Suppose that $f$ is not concave. For any $\delta>0$ we consider the lower semicontinuous function $f_\delta(x):=f(x)-\delta x^2$ with $x\in I$. It is clear that $\{f_\delta(x)\}\to f(x)$ for any $x\in I$ when $\delta\to 0$. As the concavity is preserved under pointwise convergence, there must be $\delta>0$ for which $f_\delta$ is not concave. So, we can find $x_0,y_0\in I$ with $x_0<y_0$ and $x\in(x_0,y_0)$ such that $f_\delta(x)<g(x)$, where $g:\rr\to\rr$ is the affine function with $g(x_0)=f_\delta(x_0)$ and $g(y_0)=f_\delta(y_0)$. Hence, the upper semicontinuous function $g-f_\delta:[x_0,y_0]$ attains a positive maximum at a point $m_0\in(x_0,y_0)$. For any $t\in\rr$ small enough we would have $(g-f_\delta)(m_0\pm t)\leq (g-f_\delta)(m_0)$. As a consequence
\begin{align*}
0&=g(m_0+t)+g(m_0-t)-2g(m_0)\leq f_\delta(m_0+t)+f_\delta(m_0-t)-2f_\delta(m_0)
\\
&=f(m_0+t)+f(m_0-t)-2f(m_0)-(\delta(m_0+t)^2+\delta(m_0-t)^2-2\delta m_0^2).
\end{align*}
From here we would deduce that $0<2\delta\leq\overline{D^2}f(m_0)$, which is a contradiction.
\end{proof}

The next theorem shows in particular the concavity of the function $\psi:=(I_{\Om,K})^{(n+1)/n}$ for any smooth convex domain $\Om\subset\rrn$ where anisotropic minimizers of any volume exist. Necessary conditions for the equality $\overline{D^2}\psi(v_0)=0$ are also obtained, which allows to characterize the minimizers in certain cases, see Corollary~\ref{cor:unique1} and Remark~\ref{re:unique2}.

\begin{theorem}
\label{th:fundamental}
Let $K\subset\rrn$ be a smooth strictly convex body and $\Om\subset\rrn$ a smooth convex domain. If there are anisotropic isoperimetric regions in $\Om$ of volume $v_0\in(0,V(\Om))$, then the upper second derivative of $\psi:=(I_{\Om,K})^{(n+1)/n}$ satisfies $\overline{D^2}\psi(v_0)\leq 0$. Moreover, if $\overline{D^2}\psi(v_0)=0$, and $E$ is an anisotropic minimizer of volume $v_0$, then any connected component of the regular part $\Sg$ of $\overline{\ptl E\cap\Om}$ is contained in a hyperplane or a set of the form $p_0+\la\,(\ptl K)$ for some $p_0\in\rrn$ and $\la>0$. Furthermore, the second fundamental form of $\ptl\Om$ satisfies $\emph{II}(N_K,N_K)=0$ along $\ptl\Sg$.
\end{theorem}

\begin{proof}
Consider any anisotropic minimizer $E\subset\Om$ with $V(E)=v_0$. From Proposition~\ref{prop:regularity} the interior boundary $\overline{\ptl E\cap\Om}$ is the disjoint union of a smooth hypersurface $\Sg$, possibly with boundary $\ptl\Sg=\Sg\cap\ptl\Om$, and a closed singular set $\Sg_0$ with $\mathcal{H}^{n-2}(\Sg_0)=0$. As we pointed out in Remarks~\ref{re:sentinel1} (iv), we can suppose that $E$ is an open set with $\pp_K(E,\Om)=\ar_K(\Sg)$. The equilibrium conditions in Proposition~\ref{prop:varprop} entail that the anisotropic mean curvature $H_K$ of $\Sg$ is constant and $\escpr{N_K,\xi}=0$ along $\ptl\Sg$, where $N_K$ is the outer anisotropic normal on $\Sg$ and $\xi$ is the inner unit normal on $\ptl\Om$. On the other hand, Proposition~\ref{prop:approx} provides a sequence of non-trivial functions $\{\omega_\eps\}_{\eps>0}\subset C^\infty_0(\Sg)$ such that $0\leq\omega_\eps\leq 1$, $\{\omega_\eps\}_{\eps>0}\to 1$ when $\eps\to 0$, and $\lim_{\eps\to 0}\int_\Sg|\nabla_\Sg\omega_\eps|^2\,dA=0$.

Fix $\eps>0$. Since $N_K$ is tangent to $\ptl\Om$ along $\ptl\Sg$ we can find a smooth vector field $X_\eps$ with compact support on $\rrn$, tangent to $\ptl\Om$ and satisfying $X_{\eps|\Sg}=\omega_\eps N_K$. We denote $E_t^\eps:=\phi^\eps_t(E)$ and $\Sg^\eps_t:=\phi_t^\eps(\Sg)$, where $\{\phi_t^\eps\}_{t\in\rr}$ is the one-parameter group of diffeomorphisms of $X_\eps$. For any $t\in\rr$ it is clear that $E^\eps_t\subseteq\Om$ and $\pp_K(E^\eps_t,\Om)=\ar_K(\Sg^\eps_t)$. From Corollary~\ref{cor:1stvarest}, the first variation of the functional $V_\eps(t):=V(E^\eps_t)$ is $V_\eps'(0)=\int_\Sg\omega_\eps\,\var_K\,dA\neq 0$. This allows to consider the relative anisotropic profile
\[
f_\eps(v):=\pp_K(E^\eps_{V_\eps^{-1}(v)},\Om)^{(n+1)/n}=\ar_K(\Sg^\eps_{V_\eps^{-1}(v)})^{(n+1)/n},
\]
which is locally defined around $v_0$. By having in mind the definition of $I_{\Om,K}$ and that $E$ is an anisotropic minimizer with $V(E)=v_0$, it follows that $\psi(v_0)=f_\eps(v_0)$ and $\psi\leq f_\eps$ around $v_0$. Together with the second variation formula in Corollary~\ref{cor:2ndrelprof} this implies that
\begin{equation}
\label{eq:ortiga}
\overline{D^2}\psi(v_0)\leq \overline{D^2}f_\eps(v_0)=f_\eps''(v_0)=\frac{n+1}{n}\,\ar_K(\Sg)^{1/n}\left\{\frac{nH_K^2}{\ar_K(\Sg)}+\frac{\ind_K(\omega_\eps)}{\big(\int_\Sg\omega_\eps\,\var_K\,dA\big)^2}\right\},
\end{equation}
where $\indo_K$ is the quadratic form on $C^\infty_0(\Sg)$ introduced in \eqref{eq:indexform}, i.e.
\[
\ind_K(\omega_\eps)=\int_\Sg\left\{\escpr{\nabla^K_\Sg\omega_\eps,\nabla_\Sg\omega_\eps}\,\var^2_K-\text{tr}(B_K^2)\,\omega_\eps^2\,\var_K\right\}dA-\int_{\ptl\Sg}\frac{\text{II}(N_K,N_K)}{\escpr{\nu,\xi}}\,\omega_\eps^2\,\var_K\,dL.
\]
Here $\text{II}$ is the second fundamental form of $\ptl\Om$ with respect to $\xi$, and $\nu$ is the inner conormal of $\ptl\Sg$ (recall that $\escpr{\nu,\xi}$ does not vanish along $\ptl\Sg$ by \eqref{eq:conormalest}). 

The comparison in \eqref{eq:ortiga} holds for any $\eps>0$. Hence, by taking $\limsup$ when $\eps\to 0$ and applying the dominated convergence theorem, we arrive at
\begin{equation}
\label{eq:confucio}
\overline{D^2}\psi(v_0)\leq \frac{n+1}{n}\,\ar_K(\Sg)^{1/n}\left\{\frac{nH_K^2}{\ar_K(\Sg)}+\frac{\limsup_{\eps\to 0}\ind_K(\omega_\eps)}{\ar_K(\Sg)^2}\right\}.
\end{equation}
Thus, to infer that $\overline{D^2}\psi(v_0)\leq 0$ it is enough to see that
\begin{equation}
\label{eq:desiree0}
\limsup_{\eps\to 0}\,\ind_K(\omega_\eps)\leq -nH_K^2\,\ar_K(\Sg).
\end{equation} 
For this we will analyze the different terms in $\indo_K(\omega_\eps)$.

On the one hand, by equations \eqref{eq:angrad}, \eqref{eq:endoq} and  \eqref{eq:bound}, there are constants $a,b>0$ only depending on $K$ such that
\[
a\,|\nabla_\Sg\omega_\eps|^2\leq\escpr{\nabla_\Sg^K\omega_\eps,\nabla_\Sg\omega_\eps}\leq b\,|\nabla_\Sg\omega_\eps|^2.
\]
As $\lim_{\eps\to 0}\int_\Sg|\nabla_\Sg\omega_\eps|^2\,dA=0$ and $\var_K$ is bounded, we get
\begin{equation}
\label{eq:desiree1}
\limsup_{\eps\to 0}\int_\Sg\escpr{\nabla_\Sg^K\omega_\eps,\nabla_\Sg\omega_\eps}\,\var_K^2\,dA=0.
\end{equation}
On the other hand, $\{\text{tr}(B^2_K)\,\omega_\eps^2\,\var_K\}_{\eps>0}$ is a sequence of non-negative functions (by \eqref{eq:ineq}), which pointwise converges to $\text{tr}(B^2_K)\,\var_K$ when $\eps\to 0$. From Fatou's lemma we obtain
\begin{equation}
\label{eq:desiree2}
\int_\Sg\text{tr}(B^2_K)\,\var_K\,dA\leq \liminf_{\eps\to 0}\int_\Sg\text{tr}(B^2_K)\,\omega_\eps^2\,\var_K\,dA.
\end{equation}
Next, we study the boundary term in $\indo_K(\omega_\eps)$. Observe that $\text{II}(N_K,N_K)\geq 0$ along $\ptl\Sg$ by the convexity of $\Om$. Let us check that $\escpr{\nu,\xi}>0$ along $\ptl\Sg$. Take a point $p\in\ptl\Sg$. As $\nu(p)$ is the inner conormal of $\ptl\Sg$, there is a smooth curve $\gamma:[0,\tau)\to\Sg$ such that $\gamma(0)=p$ and $\gamma'(0)=\nu(p)$. As $\Sg\subset\overline{\Om}$ and $\Om$ is convex, the function $g(t):=\escpr{\gamma(t)-p,\xi(p)}$ satisfies $g(t)\geq 0$ for any $t\in[0,\tau)$ and $g(0)=0$. This gives us $\escpr{\nu(p),\xi(p)}=g'(0)\geq 0$, which shows the claim because $\escpr{\nu(p),\xi(p)}\neq 0$. All this entails that the boundary integrand in $\indo_K(\omega_\eps)$ is nonnegative. Again from Fatou's lemma we conclude that
\begin{equation}
\label{eq:desiree3}
\int_{\ptl\Sg}\frac{\text{II}(N_K,N_K)}{\escpr{\nu,\xi}}\,\var_K\,dL\leq\liminf_{\eps\to 0}\int_{\ptl\Sg}\frac{\text{II}(N_K,N_K)}{\escpr{\nu,\xi}}\,\omega_\eps^2\,\var_K\,dL.
\end{equation}

At this point, equations \eqref{eq:desiree1}, \eqref{eq:desiree2}, \eqref{eq:desiree3} and the inequality \eqref{eq:ineq} imply that
\begin{equation}
\label{eq:desiree4}
\limsup_{\eps\to 0}\,\ind_K(\omega_\eps)\leq-\int_\Sg\text{tr}(B^2_K)\,\var_K\,dA-\int_{\ptl\Sg}\frac{\text{II}(N_K,N_K)}{\escpr{\nu,\xi}}\,\var_K\,dL\leq-nH^2_K\,\ar_K(\Sg).
\end{equation}
This shows \eqref{eq:desiree0}, so that we have proved $\overline{D^2}\psi(v_0)\leq 0$. Finally, if $\overline{D^2}\psi(v_0)=0$, by looking at \eqref{eq:confucio} and \eqref{eq:desiree4} we infer that $\text{tr}(B^2_K)=nH^2_K$ on $\Sg$, whereas $\text{II}(N_K,N_K)=0$ along $\ptl\Sg$. From Proposition~\ref{prop:umbilical} it follows that any connected component of $\Sg$ is contained in a hyperplane or a set of the form $p_0+\la\,(\ptl K)$ for some $p_0\in\rrn$ and $\la>0$. This completes the proof.
\end{proof}

Now we can prove the following result, where we do not assume additional regularity on $K$ and $\Om$.

\begin{theorem}
\label{th:concavity}
For any convex body $K\subset\rrn$ and any bounded convex domain $\Om\subset\rrn$ the function $\psi:=(I_{\Om,K})^{(n+1)/n}:[0,V(\Om)]\to\rr$ is concave.
\end{theorem}

\begin{proof}
We distinguish three cases depending on $K$ and $\Om$.

\emph{Case 1.} $K$ is a smooth strictly convex body and $\Om$ is a smooth bounded convex domain.

The function $\psi$ is continuous by Proposition~\ref{prop:conti}. Hence, to prove the claim it suffices to see that $\psi$ is concave on $(0,V(\Om))$. Since $\Om$ is bounded, Proposition~\ref{prop:existence} guarantees the existence of anisotropic minimizers for any given volume. Thus, we can apply Theorem~\ref{th:fundamental} to infer that $\overline{D^2}\psi\leq 0$ on $(0,V(\Om))$. The concavity of $\psi$ on $(0,V(\Om))$ follows from Lemma~\ref{lem:weak2nd}.

\emph{Case 2.} $K$ is a smooth strictly convex body and $\Om$ is a bounded convex domain.

We can deduce the statement from Case 1 by reasoning as Ritor\'e~\cite[Thms.~9.32, 9.33]{ritore-book}, see also Leonardi, Ritor\'e and Vernadakis~\cite[Lem.~4.1, Prop.~4.4]{unbounded}. We include below a description of the argument where we omit the details.

Consider the function $J_{\Om,K}:[0,1]\to\rr$ given by $J_{\Om,K}(t):=I_{\Om,K}(t\,V(\Om))$ for any $t\in[0,1]$. It is clear that $\psi$ is concave on $[0,V(\Om)]$ if and only if $(J_{\Om,K})^{(n+1)/n}$ is concave on $[0,1]$. To check the concavity of 
$(J_{\Om,K})^{(n+1)/n}$ we will see that it is pointwise limit of concave functions on $[0,1]$. 

From a result of Schneider~\cite{schneider-approx}, see also \cite[Sect.~5.2]{unbounded}, there is a sequence $\{\Om_i\}_{i\in\nn}$ of smooth bounded convex domains in $\rrn$ such that $\{C_i\}_{i\in\nn}$ converges in Hausdorff distance to $C$, where $C_i:=\overline{\Om}_i$ for any $i\in\nn$ and $C:=\overline{\Om}$. From Case 1 the function $(J_{\Om_i,K})^{(n+1)/n}$ is concave on $[0,1]$ for any $i\in\nn$. Next we show that the Hausdorff convergence implies that $\{J_{\Om_i,K}(t)\}\to J_{\Om,K}(t)$ for any $t\in[0,1]$ when $i\to\infty$. This is clear when $t\in\{0,1\}$. For a fixed $t\in(0,1)$ this is equivalent to that $\{I_{\Om_i,K}(v_i)\}\to I_{\Om,K}(v)$, where $v:=t\,V(\Om)$ and $v_i=t\,V(\Om_i)$ for any $i\in\nn$ (note that $\{v_i\}\to v$ when $i\to\infty$ by Lemma~\ref{lem:hausdorff} (i)). Indeed, we will prove that
\begin{equation}
\label{eq:ineqs}
I_{\Om,K}(v)\leq\liminf_{i\to\infty}I_{\Om_i,K}(v_i) \quad\text{and}\quad  I_{\Om,K}(v)\geq\limsup_{i\to\infty}I_{\Om_i,K}(v_i).
\end{equation}

By Proposition~\ref{prop:existence} we can find, for any $i\in\nn$, a set $E_i\subset\Om_i$ with $V(E_i)=v_i$ and $\pp_K(E_i,\Om_i)=I_{\Om_i,K}(v_i)$. By using the Hausdorff convergence, properties of the anisotropic perimeter, and the convexity of $\Om_i$, we infer the existence of a round open ball $B_0$ and a constant $c>0$ such that $\overline{E}_i\subset B_0$ and $\pp_K(E_i,B_0)\leq c$ for any $i\in\nn$. The compactness result in Proposition~\ref{prop:perprop} (ii) ensures that, passing to a subsequence if necessary, $\{E_i\}_{i\in\nn}$ converges to some set $E\subset B_0$ in $L^1(B_0)$. In particular $V(E)=v$ because $\{v_i\}\to v$. We can also suppose that the characteristic functions $\{\chi_{E_i}\}_{i\in\nn}$ pointwise converges to $\chi_E$ almost everywhere. Thus, we can employ Lemma~\ref{lem:hausdorff} (ii) to conclude that $E\subset\Om$ (after changing it by a set of null volume). Finally, by having in mind Lemma~\ref{lem:hausdorff} (iii), the lower semicontinuity of $\pp_K$ and that $\{E_i\}\to E$ in $L^1(B_0)$, we get
\[
I_{\Om,K}(v)\leq\pp_K(E,\Om)\leq\liminf_{i\to\infty}\pp_K(E_i,\Om_i)=\liminf_{i\to\infty}I_{\Om_i,K}(v_i),
\]
which is the first inequality in \eqref{eq:ineqs}. Next, we will obtain the second one.

After a translation (which does not change the anisotropic profiles, see Proposition~\ref{prop:perprop} (iii)) we can admit that $0\in\Om$. From Lemma~\ref{lem:hausdorff} (iv) there is $i_0\in\nn$ and a sequence $\{\la_i\}_{i\geq i_0}\subset(0,1]$ such that $\{\la_i\}\to 1$ when $i\to\infty$ and $\la_i\,C_i\subseteq C$ for any $i\geq i_0$. By using Lemma~\ref{lem:hausdorff} (ii) it is easy to see that $\{\la_i\,C_i\}_{i\in\nn}$ converges in Hausdorff distance to $C$, so that $\{V(\la_i\,\Om_i)\}\to V(\Om)$. Let $E\subset\Om$ be a set with $V(E)=v$ and $\pp_K(E,\Om)=I_{\Om,K}(v)$. As in the deformation result \cite[Thm.~1.50]{ritore-book} we can find a constant $c>0$ and a family of sets $\{E_t\}_{t\in(-\eps,\eps)}\subset\Om$ with $E_0=E$, whereas $V(E_t)=v+t$ and $\pp_K(E_t,\Om)\leq\pp_K(E,\Om)+c\,|t|=I_{\Om,K}(v)+c\,|t|$ for any $t\in(-\eps,\eps)$. Now, by taking into account that $\{V(E_t\cap\la_i\,\Om_i)\}\to V(E_t)=v+t$ and $\{\la_i^{n+1}\,v_i\}\to v$ when $i\to\infty$, we obtain the existence of $j_0\geq i_0$ and $\{t_i\}_{i\geq j_0}\subset(-\eps,\eps)$ such that $V(E_{t_i}\cap\la_i\,\Om_i)=\la_i^{n+1}\,v_i$ (this entails that $\{t_i\}\to 0$). Therefore, the set $E_i:=\la_i^{-1}\,(E_{t_i}\cap\la_i\,\Om_i)$ satisfies $E_i\subset\Om_i$ and $V(E_i)=v_i$. Moreover, the definition of $I_{\Om_i,K}$ and the second equality in Proposition~\ref{prop:perprop} (iii) give us
\begin{align*}
I_{\Om_i,K}(v_i)\leq \pp_K(E_i,\Om_i)&=\pp_K(\la_i^{-1}\,(E_{t_i}\cap\la_i\,\Om_i),\la_i^{-1}\,\la_i\,\Om_i)=\la_i^{-n}\,\pp_K(E_{t_i}\cap\la_i\,\Om_i,\la_i\,\Om_i)
\\
&\leq\la_i^{-n}\,\pp_K(E_{t_i},\Om)\leq\la_i^{-n}\,(I_{\Om,K}(v)+c\,|t_i|).
\end{align*}
From here the desired comparison in \eqref{eq:ineqs} follows by taking $\limsup$ when $i\to\infty$.

\emph{Case 3.} $K$ is a convex body and $\Om$ is a bounded convex domain.

From \cite{schneider-approx} there is a sequence of smooth strictly convex bodies $\{K_i\}_{i\in\nn}$ converging in Hausdorff distance to $K$. By Lemma~\ref{lem:hausdorff} (iv) we can find $i_0\in\nn$ and sequences $\{\theta_i\}_{i\geq i_0}\subset(0,1]$, $\{\mu_i\}_{i\geq i_0}\subset[1,\infty)$ such that $\{\theta_i\}\to 1$, $\{\mu_i\}\to 1$ and $\theta_i\,K\subseteq K_i\subseteq\mu_i\,K$ for any $i\geq i_0$. By taking support functions and having in mind equation~\eqref{eq:hK} we get $\theta_i\,h_{K}\leq h_{K_i}\leq\mu_i\,h_{K}$ for any $i\geq i_0$. Equation \eqref{eq:perfinite} leads to
\[
\theta_i\,\pp_K(E,\Om)\leq\pp_{K_i}(E,\Om)\leq\mu_i\,\pp_K(E,\Om)
\]
for any set $E$ of locally finite perimeter in $\Om$. As a consequence, for any $i\geq i_0$ we have
\[
\theta_i\,I_{\Om,K}\leq I_{\Om,K_i}\leq\mu_i\,I_{\Om,K},
\]
so that $\{I_{\Om,K_i}\}_{i\in\nn}$ pointwise converges to $I_{\Om,K}$ in $[0,V(\Om)]$. Therefore the concavity of $(I_{\Om,K_i})^{(n+1)/n}$, deduced from Case 2, allows to conclude that $(I_{\Om,K})^{(n+1)/n}$ is concave as well.
\end{proof}

From the previous theorem we can easily derive some analytical properties of $I_{\Om,K}$. 

\begin{corollary}
\label{cor:regprop}
Let $K\subset\rrn$ be a convex body and $\Om\subset\rrn$ a bounded convex domain.
\begin{itemize}
\item[(i)] The anisotropic isoperimetric profile $I_{\Om,K}:[0,V(\Om)]\to\rr$ is a concave function.
\item[(ii)] $I_{\Om,K}$ is strict subbaditive on $(0,V(\Om))$, i.e., $I_{\Om,K}(v_1+v_2)<I_{\Om,K}(v_1)+I_{\Om,K}(v_2)$ for any $v_1,v_2\in(0,V(\Om))$ such that $v_1+v_2\in(0,V(\Om))$.
\item[(iii)] The left and right derivatives $(I_{\Om,K})'_-$ and $(I_{\Om,K})'_+$ exist and are non-increasing on $(0,V(\Om))$. Moreover, $I_{\Om,K}$ is absolutely continuous and differentiable except on an at most countable set.
\item[(iv)] Suppose that $K$ is smooth and strictly convex. Then, for any $v_0\in(0,V(\Om))$ and any anisotropic minimizer $E\subset\Om$ with $V(E)=v_0$ we have
\[
(I_{\Om,K})'_-(v_0)\geq-nH_K\geq (I_{\Om,K})'_+(v_0),
\]
where $H_K$ is the associated anisotropic mean curvature of $E$, see Remarks~\ref{re:sentinel1} (iii) and Remark~\ref{re:sentinel2}. In particular, if $I_{\Om,K}$ is differentiable at $v_0$, then $H_K$ does not depend on $E$. 
\item[(v)] If $K$ is also centrally symmetric about $0$, then $I_{\Om,K}$ is non-decreasing on $[0,V(\Om)/2]$. Furthermore, for any anisotropic minimizer $E$ in $\Om$ with $V(E)<V(\Om/2)$, we have $H_K\leq 0$.
\end{itemize}  
\end{corollary}

\begin{proof}
The function $\psi:=(I_{\Om,K})^{(n+1)/n}$ is concave by Theorem~\ref{th:concavity}. This, together with the fact that $\psi(0)=0$, entails that $\psi$ is subadditive on $[0,V(\Om)]$. On the other hand, the function $\var:[0,\infty)\to\rr$ given by $\var(x):=x^{n/(n+1)}$ is increasing and strictly concave. The statements (i) and (ii) come because $I_{\Om,K}=\var\circ\psi$ and $\var$ is strictly subadditive on $(0,\infty)$. The regularity of $I_{\Om,K}$ stated in (iii) is consequence of its concavity, see for instance \cite[Ch.~1, \S 4]{bourbaki} and \cite[Cor.~24.2.1]{rockafellar}. 

Let us prove (iv). Denote by $\Sg$ the regular part of $\ptl E\cap\Om$, where $E\subset\Om$ is any minimizer with $V(E)=v_0$. We repeat the construction in the proof of Theorem~\ref{th:fundamental} with any $\omega\in C^\infty_0(\Sg)$ such that $\int_\Sg\omega\,\var_K\,dA\neq 0$. If $g(v):=\ar_K(V^{-1}(v))$ stands for the relative profile of the associated deformation of $E$, then $I_{\Om,K}(v)\leq g(v)$ for $v$ close enough to $v_0$ and $I_{\Om,K}(v_0)=g(v_0)$. From here we obtain
\[
(I_{\Om,K})'_-(v_0)\geq g'(v_0)=\frac{\ar_K'(0)}{V'(0)}\geq(I_{\Om,K})'_+(v_0),
\]
and the desired conclusion follows from Corollary~\ref{cor:1stvarest}. Finally, when $K$ is centrally symmetric we know from Remark~\ref{re:centrally} that $I_{\Om,K}$ is symmetric about $V(\Om)/2$. This gives us
\[
(I_{\Om,K})'_-(V(\Om)/2)=-(I_{\Om,K})'_+(V(\Om)/2)\geq-(I_{\Om,K})'_-(V(\Om)/2),
\]
so that $(I_{\Om,K})'_-\geq(I_{\Om,K})'_-(V(\Om)/2)\geq 0$ on $(0,V(\Om)/2]$. Since $I_{\Om,K}$ is absolutely continuous this implies that $I_{\Om,K}$ is non-decreasing on $[0,V(\Om)/2]$. So, if $E$ is any anisotropic minimizer with $v_0=V(E)<V(\Om)/2$, then $-nH_K\geq(I_{\Om,K})'_+(v_0)\geq 0$.
\end{proof}

A direct application of the strict subadditivity of the anisotropic profile is the connectivity of the anisotropic minimizers.

\begin{theorem}
\label{th:connectset}
For any convex body $K\subset\rrn$ and any bounded convex domain $\Om\subset\rrn$ the anisotropic isoperimetric regions in $\Om$ are connected.
\end{theorem}

\begin{proof}
Let $E\subset\Om$ be an anisotropic minimizer with $V(E)=v_0$. Suppose that $E$ is not connected and take any connected component $E_1$ of $E$. If we denote $v_1:=V(E_1)$ then $v_1<v_0$ and we have
\begin{align*}
\pp_K(E,\Om)&=I_{\Om,K}(v_0)=I_{\Om,K}(v_1+v_0-v_1)<I_{\Om,K}(v_1)+I_{\Om,K}(v_0-v_1)
\\
&\leq\pp_K(E_1,\Om)+\pp_K(E\setminus E_1,\Om)=\pp_K(E,\Om),
\end{align*}
where the strict inequality comes from Corollary~\ref{cor:regprop} (ii). This produces the desired contradiction.
\end{proof}

We now study the connectivity of the regular part in the interior boundary of an anisotropic minimizer. Observe that this is a stronger property than the connectivity of the minimizers. Our analysis relies on the second variation formula in Proposition~\ref{prop:2ndareavol}, so that additional regularity conditions on $K$ and $\Om$ are required.

\begin{theorem}
\label{th:connectbound}
Let $K\subset\rrn$ be a smooth strictly convex body and $\Om\subset\rrn$ a smooth convex domain. Consider an anisotropic minimizer $E\subset\Om$ and denote by $\Sg$ the regular part of $\overline{\ptl E\cap\Om}$. Then, either $\Sg$ is connected, or any connected component of $\Sg$ is contained in a hyperplane and satisfies equality $\emph{II}(N_K,N_K)=0$ along the points in $\ptl\Sg$.
\end{theorem}

\begin{proof}
From Proposition~\ref{prop:varprop} the anisotropic mean curvature $H_K$ is constant on $\Sg$ and $\escpr{N_K,\xi}=0$ along $\ptl\Sg$. Suppose that $\Sg$ is disconnected and take any two components $\Sg_1\neq\Sg_2$. To prove the claim we will find a set of the same volume as $E$ and strictly less anisotropic perimeter unless $\Sg$ verifies the announced conditions.

Choose a non-trivial function $\omega_i\in C^\infty_0(\Sg_i)$ with $\omega_i\geq 0$, and extend it by $0$ over $\Sg\setminus\Sg_i$. Let $X_i$ be a smooth vector field with compact support on $\rrn$, tangent to $\ptl\Om$, and satisfying $X_{i|\Sg}=\omega_i N_K$. We can construct $X_i$ in such a way that $\text{supp}(X_i)\subset U_i$ and $U_i\cap (\Sg\setminus\Sg_i)=\emptyset$ for some bounded open set $U_i\subset\rrn$. Moreover,  as $\Sg_1\cap\Sg_2=\emptyset$, we can suppose that $U_1\cap U_2=\emptyset$. We denote by $\{(\phi_1)_t\}_{t\in\rr}$ and $\{(\phi_2)_s\}_{s\in\rr}$ the one-parameter groups of diffeomorphisms associated to $X_1$ and $X_2$, respectively. The corresponding area and volume functionals are given by $(\ar_K)_1(t):=\ar_K((\phi_1)_t(\Sg))$, $(\ar_K)_2(s):=\ar_K((\phi_2)_s(\Sg))$, $V_1(t):=V((\phi_1)_t(E))$ and $V_2(s):=V((\phi_2)_s(E))$, for any $t,s\in\rr$.

Consider the set $E_{t,s}:=((\phi_2)_s\circ(\phi_1)_t)(E)$ and the function $V(t,s):=V(E_{t,s})$. It is clear that $(\ptl V/\ptl t)(t,0)=V_1'(t)$, $(\ptl V/\ptl s)(0,s)=V_2'(s)$, $(\ptl^2 V/\ptl t^2)(0,0)=V_1''(0)$ and $(\ptl^2 V/\ptl s^2)(0,0)=V_2''(0)$. We also have $(\ptl^2 V/\ptl s\ptl t)(0,0)=0$ because $\text{supp}(X_i)\subset U_i$ and $U_1\cap U_2=\emptyset$. From the first variation of volume in Corollary~\ref{cor:1stvarest}, we get
\[
\frac{\ptl V}{\ptl s}(0,0)=V_2'(0)=\int_\Sg\omega_2\,\var_K\,dA>0.
\] 
Thus, the implicit function theorem entails the existence of a function $s(t)$ with $t\in(-\delta,\delta)$ such that $s(0)=0$ and $V(t,s(t))=V(E)$ for any $t\in(-\delta,\delta)$. By differentiating with respect to $t$, we obtain
\[
s'(0)=-\frac{V_1'(0)}{V_2'(0)}\quad\text{and}\quad s''(0)=-\frac{V_1''(0)}{V_2'(0)}-\frac{V_1'(0)^2\,V_2''(0)}{V_2'(0)^3}.
\]

Define $E_t:=E_{t,s(t)}$ for any $t\in(-\delta,\delta)$. This a set in $\Om$ having the same regularity properties as $E$ and the same volume. In particular $\pp_K(E_t,\Om)=\ar_K(\Sg_t)$, where $\Sg_t:=((\phi_2)_{s(t)}\circ(\phi_1)_t)(\Sg)$. By using again that $\text{supp}(X_i)\subset U_i$ and $U_1\cap U_2=\emptyset$ it follows that, for any $t$ small enough, the function $\ar_K(t):=\ar_K(\Sg_t)$ satisfies
\[
\ar_K(t)=(\ar_K)_1(t)+(\ar_K)_2(s(t))+c,
\]
for some $c\in\rr$. The fact that $E$ is an anisotropic isoperimetric region implies that $\ar_K(t)$ attains its minimum at $t=0$, so that $\ar_K'(0)=0$ and $\ar_K''(0)\geq 0$. By differentiating twice into the previous equation, substituting the values of $s'(0)$ and $s''(0)$, and having in mind Corollary~\ref{cor:1stvarest}, we obtain
\begin{align*}
\ar_K''(0)&=(\ar_K)_1''(0)+s''(0)\,(\ar_K)_2'(0)+s'(0)^2\,(\ar_K)_2''(0)
\\
&=((\ar_K)_1+nH_KV_1)''(0)+\frac{V_1'(0)^2}{V_2'(0)^2}\,((\ar_K)_2+nH_KV_2)''(0).
\end{align*}
Thus, the second variation formula in Proposition~\ref{prop:2ndareavol} leads to
\begin{equation}
\label{eq:goty24}
0\leq\ar_K''(0)=\indo_K(\omega_1)+\frac{(\int_{\Sg}\omega_1\,\var_K\,dA)^2}{(\int_\Sg\omega_2\,\var_K\,dA)^2}\,\indo_K(\omega_2),
\end{equation}
where $\indo_K$ is the quadratic form introduced in \eqref{eq:indexform}. The previous inequality is then valid for any non-trivial functions $\omega_i\in C^\infty_0(\Sg_i)$ with $\omega_i\geq 0$. 

Now, consider a sequence $\{\omega_\eps\}_{\eps>0}$ as in Proposition~\ref{prop:approx} and set $\omega_i:=\omega_{\eps|\Sg_i}$ for any $i=1,2$. By taking $\limsup$ when $\eps\to 0$ in \eqref{eq:goty24} and reasoning as in the proof of \eqref{eq:desiree4}, we deduce 
\begin{align*}
0&\leq-\int_{\Sg_1}\text{tr}(B^2_K)\,\var_K\,dA-\int_{\ptl\Sg_1}\frac{\text{II}(N_K,N_K)}{\escpr{\nu,\xi}}\,\var_K\,dL
\\
&-\frac{\ar_K(\Sg_1)^2}{\ar_K(\Sg_2)^2}\,\left(\int_{\Sg_2}\text{tr}(B^2_K)\,\var_K\,dA+\int_{\ptl\Sg_2}\frac{\text{II}(N_K,N_K)}{\escpr{\nu,\xi}}\,\var_K\,dL\right).
\end{align*}
From the convexity of $\Om$ we know that $\text{II}(N_K,N_K)\geq 0$ and $\escpr{\nu,\xi}>0$ (this was shown just before equation~\eqref{eq:desiree3}). We conclude from above that $\text{tr}(B^2_K)=0$ on $\Sg_i$ and $\text{II}(N_K,N_K)=0$ along $\ptl\Sg_i$ for any $i=1,2$. The proof finishes by invoking Proposition~\ref{prop:umbilical}.
\end{proof}

\begin{remark}
The minimizing property of $E$ is only employed to ensure that $\ar_K'(0)=0$ and $\ar_K''(0)\geq 0$. So, the theorem also holds for a set $E$ of finite volume and anisotropic perimeter, which is a second order minima of $\pp_K(\cdot\,,\Om)$ under a volume constraint, satisfies the regularity properties in Proposition~\ref{prop:regularity} (it is enough that the singular set $\Sg_0$ is relatively closed with $\ar(\Sg_0)=0$) and the approximation result in Proposition~\ref{prop:approx}.
\end{remark}

Theorem~\ref{th:connectbound} does not prevent connected minimizers with $\Sg$  disconnected. For instance, in a strip $\Om\subset\rr^2$, the minimizers of the isotropic length for big enough area are separated by two segments orthogonal to $\ptl\Om$. The next result shows additional conditions ensuring that $\Sg$ is connected.

\begin{corollary}
\label{cor:connected24}
Let $K\subset\rrn$ be a smooth strictly convex body and $\Om\subset\rrn$ a smooth convex domain. If $\Om$ is strictly convex, or $K$ is centrally symmetric about $0$ and $\Om$ is bounded, then the regular part $\Sg$ in the interior boundary $\overline{\ptl E\cap\Om}$ of any anisotropic minimizer $E\subset\Om$ must be connected.
\end{corollary}

\begin{proof}
Suppose that $\Sg$ is disconnected. By Theorem~\ref{th:connectbound} any component $\Sg'\subset\Sg$ is contained in a hyperplane $P'$ and satisfies $\text{II}(N_K,N_K)=0$ along the points in $\ptl\Sg$. Note that the singular set $\Sg_0:=\overline{\ptl E\cap\Om}\setminus\Sg$ is contained in $\ptl\Om$. Otherwise, we would find an $(n-1)$-dimensional submanifold $S$ (contained in the intersection of two different hyperplanes) such that $S\subseteq\Sg_0$, which contradicts the equality $\mathcal{H}^{n-2}(\Sg_0)=0$ in Proposition~\ref{prop:regularity}. The fact that $\Sg_0\subseteq\ptl\Om$ entails that $\overline{\Sg'}=P'\cap\overline{\Om}$. Observe that $P'\cap\ptl\Om\neq\emptyset$ because $\ar_K(\Sg')<\infty$. As $\Om$ is convex, $P'$ meets $\ptl\Om$ transversally and so, $P'\cap\ptl\Om$ is an $(n-1)$-dimensional submanifold. By using again that $\mathcal{H}^{n-2}(\Sg_0)=0$ we infer that $(P'\cap\ptl\Om)\setminus\Sg_0\neq\emptyset$. This implies that $\Sg'\cap\ptl\Om\neq\emptyset$, so that $\ptl\Sg'\neq\emptyset$. Thus, from the equality $\text{II}(N_K,N_K)=0$ we deduce that $\Om$ cannot be strictly convex. Neither is it possible that $K$ is centrally symmetric about $0$ and $\Om$ is bounded. Otherwise, by Remark~\ref{re:centrally} and Theorem~\ref{th:connectset} we would deduce that $E$ and $\Om\setminus\overline{E}$ are both connected. As any hyperplane $P'$  separates $\Om$ into two connected components, this would enforce $\Sg$ to be connected, a contradiction.
\end{proof}

\begin{remark}
The aforementioned example of a planar strip illustrates that the corollary need not hold when $K$ is centrally symmetric about $0$ but $\Om$ is unbounded. In relation to this, Sternberg and Zumbrun~\cite[Thm.~2.6]{sz} proved that if $E$ is a bounded local minimizer of the isotropic perimeter for fixed volume in a smooth convex domain $\Om$, then $\overline{\Sg}$ is connected or consists of components inside parallel hyperplanes with that part of $\Om$ between any two such components consisting of a cylinder.
\end{remark}

Our next aim is to establish a sharp anisotropic isoperimetric inequality for a convex domain $\Om$. More precisely, we will obtain the profile comparison $I_{\Om,K}\leq I_{\cc,K}$, where $\cc$ is a certain open convex cone. Such an estimate is suggested by the fact that $(I_{\Om,K})^{(n+1)/n}$ is concave (Theorem~\ref{th:concavity}) whereas $(I_{\cc,K})^{(n+1)/n}$ is linear (equation~\eqref{eq:isoprofcones}). Based on this, we only have to find a cone $\cc$ so that the desired comparison holds for small volumes. This is done in the next proposition for open sets satisfying a mild convexity condition. We will denote $K_{p_0,\la}:=p_0+\la K$ for any $p_0\in\rrn$ and $\la>0$.

\begin{proposition}
\label{prop:1stcomp}
Let $K\subset\rrn$ be a convex body and $\Om\subset\rrn$ an open set having a local supporting half-space at a point $p_0\in\ptl\Om$. Then, for any open convex cone $\cc\subset\rrn$ with vertex $p_0$ and containing $\Om$, there is $v_0>0$ such that
\[
I_{\Om,K}(v)\leq I_{\cc,K}(v), \quad\text{for any } v\in (0,v_0).
\]
Moreover, if equality holds for some $v\in(0,v_0)$ and we take $\la>0$ for which $V(K_{p_0,\la}\cap\Om)=v$, then $K_{p_0,\la}\cap\Om$ is an anisotropic minimizer in $\Om$ and $K_{p_0,\la}\cap\Om=K_{p_0,\la}\cap\cc$.
\end{proposition}

\begin{proof}
To simplify the notation we set $V_A(\la):=V(K_{p_0,\la}\cap A)$ and $\pp_A(\la):=\pp_K(K_{p_0,\la},A)$ for any open set $A\subseteq\rrn$ and $\la>0$. Consider the open cone $\cc_\la$ defined as the union of all open half-lines joining $p_0$ with the Lipschitz hypersurface $\ptl K_{p_0,\la}\cap\Om$. The local convexity of $\Om$ at $p_0$ implies the existence of $\la_0>0$ such that $K_{p_0,\la}\cap\cc_\la\subseteq K_{p_0,\la}\cap\Om$ for any $\la\in(0,\la_0)$. It is also clear that $K_{p_0,\la}\cap\Om\subseteq K_{p_0,\la}\cap\cc$ for any $\la>0$. Thus, we have
\[ 
V_{\cc_\la}(\la)\leq V_\Om(\la)\leq V_\cc(\la),\quad\text{for any }\la\in(0,\la_0).
\] 
On the other hand, from equation~\eqref{eq:relationcone} and Proposition~\ref{prop:perprop} (iii) we get
\[
\frac{\pp_{\cc_\la}(\la)}{V_{\cc_\la}(\la)}=\frac{n+1}{\la}=\frac{\pp_\cc(\la)}{V_C(\la)}, \quad\text{for any }\la\in(0,\la_0).
\]
By using again Proposition~\ref{prop:perprop} (iii) and the isoperimetric property of the sets $\la K\cap (\cc-p_0)$ in the convex cone $\cc-p_0$, we infer that $K_{p_0,\la}\cap\cc$ is an anisotropic minimizer in $\cc$ for any $\la>0$. All this together with the definition of $I_{\Om,K}$ and equation~\eqref{eq:isocones} entails
\begin{align*}
I_{\Om,K}(V_\Om(\la))&\leq\pp_\Om(\la)=\pp_{\cc_\la}(\la)=\frac{n+1}{\la}\,V_{\cc_\la}(\la)\leq\frac{\pp_\cc(\la)}{V_\cc(\la)}\,V_\Om(\la)
\\
&=\frac{c(\cc,K)}{V_\cc(\la)^{1/(n+1)}}\,V_\Om(\la)\leq c(\cc,K)\,V_\Om(\la)^{n/(n+1)}=I_{\cc,K}(V_\Om(\la)),
\end{align*}
for any $\la\in(0,\la_0)$. This provides the desired comparison because $V_\Om(\la)$ is a continuous non-decreasing function. Finally, if equality holds for some $v=V_\Om(\la)$, then $I_{\Om,K}(v)=\pp_\Om(\la)$ and $V_{\cc_\la}(\la)=V_\Om(\la)=V_\cc(\la)$. Thus, $K_{p_0,\la}\cap\Om$ is isoperimetric in $\Om$ and $(K_{p_0,\la}\cap\cc)\setminus (K_{p_0,\la}\cap\Om)$ has null volume. By convexity this yields $K_{p_0,\la}\cap\Om=K_{p_0,\la}\cap\cc$, which proves the claim.
\end{proof}

The previous statement holds in particular when $\Om$ is convex. In this case the concavity of $(I_{\Om,K})^{(n+1)/n}$ leads to a global profile comparison.

\begin{theorem}
\label{th:maincomp}
Let $K\subset\rrn$ be a convex body and $\Om\subset\rrn$ a bounded convex domain. Then, for any open convex cone $\cc\subset\rrn$ containing $\Om$ and having vertex at a point $p_0\in\ptl\Om$, we have
\[
I_{\Om,K}(v)\leq I_{\cc,K}(v), \quad\text{for any } v\in [0,V(\Om)].
\]
Suppose that equality holds for $v_0\in(0,V(\Om))$ and take $\la_0>0$ for which $V(K_{p_0,\la_0}\cap\Om)=v_0$. Then, for any $\la\in(0,\la_0]$, the set $K_{p_0,\la}\cap\Om$ is an anisotropic minimizer in $\Om$ and coincides with $K_{p_0,\la}\cap\cc$. As a consequence, there is an open set $U\subset\rrn$ with $p_0\in U$ and $U\cap\ptl\Om=U\cap\ptl\cc$.
\end{theorem}

\begin{proof}
We can apply Proposition~\ref{prop:1stcomp} in order to find  $v_0\in(0,V(\Om))$ such that
\[
(I_{\Om,K})^{(n+1)/n}\leq (I_{\cc,K})^{(n+1)/n} \quad\text{in } [0,v_0]. 
\]
The function at the left is concave by Theorem~\ref{th:concavity}, whereas the one at the right is linear by equation~\eqref{eq:isoprofcones}. Thus, the previous inequality holds in $[0,V(\Om)]$. In case of equality for some $v_0\in(0,V(\Om))$, the concavity of $(I_{\Om,K})^{(n+1)/n}-(I_{\cc,K})^{(n+1)/n}$ implies $I_{\Om,K}\geq I_{\cc,K}$ in $[0,v_0]$. So, $I_{\Om,K}=I_{\cc,K}$ in $[0,v_0]$, and the conclusion follows from the discussion of equality in Proposition~\ref{prop:1stcomp}.
\end{proof}

\begin{example}
The result applies when $\cc$ is the \emph{support cone} of $\Om$ at a point $p\in\ptl\Om$. This is the smallest cone $\cc_p$ with vertex at $p$ and containing $\Om$; it is given by the union of all open half-lines joining $p$ with the points of $\Om$. Since $\cc_p$ is open and convex for any $p\in\ptl\Om$, we obtain
\[
I_{\Om,K}(v)\leq I_{\cc_p,K}(v)=(n+1)\,V(K\cap(\cc_p-p))^{1/(n+1)}\,v^{n/(n+1)}, \quad\text{for any } v\in[0,V(\Om)],
\]
where we have used that $I_{\cc_p,K}=I_{\cc_p-p,K}$ and equation~\eqref{eq:isoprofcones}. On the other hand, we can proceed as in the isotropic case, see Ritor\'e and Vernadakis~\cite[Lem.~6.1]{unbounded}, to show that the function
\[
\theta(p):=V(K\cap(\overline{\cc}_p-p))=V((p+K)\cap\overline{\cc}_p)
\]
is lower semicontinuous along $\ptl\Om$. Hence, its minimum value $\theta_0$ is attained on $\ptl\Om$, and we deduce
\[
I_{\Om,K}(v)\leq (n+1)\,\theta_0^{1/(n+1)}\,v^{n/(n+1)}, \quad\text{for any } v\in[0,V(\Om)].
\]
Moreover, if equality holds for some $v_0\in (0,V(\Om))$ and we take any point $p_0\in\ptl\Om$ where $\theta(p)=\theta_0$, then there is an open set $U\subset\rrn$ with $p_0\in U$ and $U\cap\ptl\Om=U\cap\ptl\cc_{p_0}$.
\end{example}

\begin{example}
\label{ex:halfspace}
For any open supporting half-space $\hh$ of $\Om$ at a point $p_0\in\ptl\Om$, the theorem leads to the comparison $I_{\Om,K}\leq I_{\hh,K}$ on $[0,V(\Om)]$. Moreover, if equality holds for some $v_0\in(0,V(\Om))$, then there is an open neighborhood of $p_0$ in $\ptl\Om$ which is part of a hyperplane. This allows us to deduce the strict inequality $I_{\Om,K}<I_{\hh,K}$ in $(0,V(\Om))$ when $\Om$ is strictly convex.

As we observed in Remark~\ref{re:halfspace} the isoperimetric profile $I_{\hh,K}$ may depend on the half-space $\hh$. By equation~\eqref{eq:isoprofcones} such profile is determined by $V(K\cap(\hh-p))$. Note that $\hh-p$ is a half-space of the form $\hh_w:=\{p\in\rrn\,;\,\escpr{p,w}>0\}$ for some $w\in\sph^n$. As the function $w\mapsto V(K\cap\hh_w)$ attains its minimum at $w_0\in\sph^n$ then, by taking any $p_0\in\ptl\Om$ where an open supporting half-space of $\Om$ coincides with $\hh_{w_0}$, we infer that
$I_{\Om,K}\leq I_{\hh_{w_0},K}$ in $[0,V(\Om)]$. 

In the case where $K$ is centrally symmetric about $0$ we have that $I_{\hh,K}$ does not depend on the half-space $\hh$. In this situation we get the strict inequality $I_{\Om,K}<I_{\hh,K}$ on $(0,V(\Om))$. This is because, in case of equality, we would deduce from Theorem~\ref{th:maincomp} that any $p\in\ptl\Om$ has a neighborhood in $\ptl\Om$ which is part of a hyperplane, thus contradicting the compactness of $\overline{\Om}$.
\end{example}

To finish this work we would like to indicate how some of our results can be combined to characterize the anisotropic isoperimetric regions in $\rrn$ when $K$ is a smooth strictly convex body. This provides a different proof of the known fact that, up to translations, dilations about $0$, and sets of volume zero, a convex body $K$ uniquely minimizes the perimeter $\pp_K$ in $\rrn$ among sets of the same volume, see Taylor~\cite{taylor-unique,taylor}, Fonseca and Müller~\cite{fonseca-muller}, or Brothers and Morgan~\cite{brothers-morgan}. 

\begin{corollary}
\label{cor:unique1}
Let $K\subset\rrn$ be a smooth strictly convex body. If $E$ is any anisotropic minimizer in $\rrn$, then $\ptl E=p_0+\la(\ptl K)$ for some $p_0\in\rrn$ and $\la>0$.
\end{corollary}

\begin{proof}
From equation~\eqref{eq:isocones}, which is a consequence of Proposition~\ref{prop:perprop} (iii), we see that $\psi:=(I_{\rrn,K})^{(n+1)/n}$ is an increasing linear function. Thus $\overline{D^2}\psi=0$ and Theorem~\ref{th:fundamental} implies that any connected component of the regular part $\Sg$ of $\ptl E$ is contained in a hyperplane or in $p_0+\la\,(\ptl K)$ for some $p_0\in\rrn$ and $\la>0$. By Corollary~\ref{cor:regprop} (iv) the outer anisotropic mean curvature of $\Sg$ satisfies $H_K<0$, which prevents connected components inside hyperplanes. Hence, Theorem~\ref{th:connectbound} entails that $\Sg$ is a connected hypersurface with $\Sg\subseteq p_0+\la\,(\ptl K)$. As the unit normal $N$ to $\Sg$ can be continuously extended to $\overline{\Sg}=\ptl E$, we infer from \cite[Thm.~4.11]{giusti} that the singular set $\Sg_0\subset\ptl E$ is empty. So, $\Sg=\ptl E$ and the inclusion $\Sg\subseteq p_0+\la\,(\ptl K)$ allows to conclude that $\ptl E$ is compact and $\ptl E=p_0+\la\,(\ptl K)$. This proves the claim.
\end{proof}

\begin{remark}
\label{re:unique2}
Similarly, the interior boundary of any anisotropic isoperimetric region within an open half-space $\hh\subset\rrn$ with $0\in\ptl\hh$ has the form $p_0+\la(\ptl K\cap\overline{\hh})$ for some $p_0\in\ptl\hh$. Indeed, the same arguments apply to any open convex cone $\cc$ with vertex at $0$ and such that $\ptl\cc\setminus\{0\}$ is smooth. This is a different approach to a more general result of Dipierro, Poggesi and Valdinoci \cite[Thm.~4.2]{dipierro-poggesi-valdinoci}, who established the uniqueness of $K\cap\cc$ as an anisotropic minimizer in $\cc$ by means of the optimal transport technique employed by Figalli and Indrei in the isotropic case \cite[Thm.~2.2]{figalli-indrei}. 
\end{remark}

\appendix
\section{Anisotropic minimal hypersurfaces with free boundary}
\label{app:minimal}
\setcounter{theorem}{0}
\setcounter{equation}{0}  
\setcounter{subsection}{0}
\noindent

The variational formulas in Section~\ref{sec:2ndvar} allow also to study the critical points and the second order minima of the anisotropic area (without a volume constraint) inside Euclidean domains. We include here some useful facts about this topic.

Let $K$ be a smooth strictly convex body and $\Om\subset\rrn$ a smooth open set. Take a two-sided smooth hypersurface $\Sg\subset\overline{\Om}$ with non-empty boundary $\ptl\Sg=\Sg\cap\ptl\Om$. We consider the variation $\{\Sg_t\}_{t\in\rr}$ of $\Sg$ associated to a complete smooth vector field $X$ on $\rrn$ which is tangent to $\ptl\Om$. Observe that $\ptl\Sg_t=\Sg_t\cap\ptl\Om$ for any $t\in\rr$, so that $\ptl\Sg$ does not leave $\ptl\Om$ along the variation. When $X$ has compact support on $\Sg$ we can define the area functional $\ar_K(t)$ as in \eqref{eq:areafunct}. Recall that this requires fixing a smooth unit normal vector field $N$ on $\Sg$.

In these conditions, $\Sg$ is \emph{anisotropic stationary} if $\ar_K'(0)=0$ for any \emph{admissible} vector field, i.e., any complete smooth vector field $X$ which is tangent to $\ptl\Om$ and has compact support on $\Sg$. From the first variation formula in \eqref{eq:1stareagen}, we get that $\Sg$ is anisotropic stationary if and only if $H_K=0$ on $\Sg$ and $\escpr{N_K,\xi}=0$ along $\ptl\Sg$. By following classical terminology we will say that $\Sg$ is an \emph{anisotropic minimal hypersurface with free boundary in $\ptl\Om$}. Such a hypersurface is \emph{stable} if $\ar_K''(0)\geq 0$ for any admissible vector field. Given a function $\omega\in C^\infty_0(\Sg)$, the equality $\escpr{N_K,\xi}=0$ along $\ptl\Sg$ allows to find an admissible vector field $X$ with $X_{|\Sg}=\omega N_K$. Then, we can proceed as in the proof of Proposition~\ref{prop:2ndareavol} to deduce that $\ar_K''(0)=\ind_K(\omega)$, where $\ind_K$ is the quadratic form in \eqref{eq:indexform}. Therefore, if $\Sg$ is a stable anisotropic minimal surface with free boundary, then $\ind_K(\omega)\geq 0$ for any $\omega\in C^\infty_0(\Sg)$. 

The boundary term in the expression of $\ind_K(\omega)$ suggests that the stability inequality $\indo_K(\omega)\geq 0$ is more restrictive when $\Om$ is convex. In this direction, we can obtain the following statement extending a well-known result in the isotropic setting.

\begin{theorem}
\label{th:minimal}
Let $K$ be a smooth strictly convex body and $\Om\subset\rrn$ a smooth convex domain. Suppose that $\Sg\subset\overline{\Om}$ is two-sided smooth hypersurface with non-empty boundary $\ptl\Sg=\Sg\cap\ptl\Om$. If $\Sg$ is a stable and parabolic anisotropic  minimal hypersurface with free boundary in $\ptl\Om$, then $\Sg$ is contained in a hyperplane and $\emph{II}(N_K,N_K)=0$ along $\ptl\Sg$.
\end{theorem}

\begin{proof}
The parabolicity of $\Sg$ is characterized by the fact that there is a sequence $\{\omega_\eps\}_{\eps>0}\sub C^\infty_0(\Sg)$ with the same properties as in Proposition~\ref{prop:approx}, see for instance \cite[Prop.~4.1]{troyanov}. From the stability of $\Sg$ we have $\ind_K(\omega_\eps)\geq 0$ for any $\eps>0$. By taking $\limsup$ when $\eps\to 0$ and reasoning as in the proof of Theorem~\ref{th:fundamental}, we arrive at
\[
0\leq \limsup_{\eps\to 0}\,\ind_K(\omega_\eps)\leq-\int_\Sg\text{tr}(B^2_K)\,\var_K\,dA-\int_{\ptl\Sg}\frac{\text{II}(N_K,N_K)}{\escpr{\nu,\xi}}\,\var_K\,dL.
\]
From here the claim follows by having in mind Proposition~\ref{prop:umbilical} and the convexity of $\Om$.
\end{proof}

\begin{remark}
The parabolicity property is satisfied when $\Sg$ is compact, or it is complete with $\ar_K(\Sg)<\infty$. In particular, any compact anisotropic minimal hypersurface with free boundary in $\ptl\Om$ is unstable if $\Om$ is strictly convex. This need not hold for a general convex domain $\Om$. For instance, take a solid cylinder $\Om:=D\times\rr$ over a smooth convex body $D\subset\rr^n$, and let $\Sg$ be a non-vertical hyperplane intersected with $\overline{\Om}$. This is a compact anisotropic minimal hypersurface. Moreover, $\Sg$ has free boundary in $\ptl\Om$ if and only if the convex body $K$ is chosen so that the corresponding vector $N_K$ is vertical. Since $\Sg$ is inside a hyperplane and $\text{II}(N_K,N_K)=0$ along $\ptl\Sg$, we infer from~\eqref{eq:indexform} that $\indo_K(\omega)=0$ for any $\omega\in C^\infty_0(\Sg)$. This indicates that $\Sg$ should be stable. Indeed, a calibration argument as in \cite[Sect.~4]{morgan-cone-regular}, see also \cite[p.~111]{rosales-bernstein}, entails that $\Sg$ minimizes $\ar_K$ among all compact hypersurfaces $\Sg'\subset\overline{\Om}$ with $\ptl\Sg'=\Sg'\cap\ptl\Om$.
\end{remark}

\providecommand{\bysame}{\leavevmode\hbox to3em{\hrulefill}\thinspace}
\providecommand{\MR}{\relax\ifhmode\unskip\space\fi MR }
\providecommand{\MRhref}[2]{%
  \href{http://www.ams.org/mathscinet-getitem?mr=#1}{#2}
}
\providecommand{\href}[2]{#2}

\end{document}